# Logics preserving degrees of truth from varieties of residuated lattices

Bou, F., Esteva, F., Font, J. M., Gil, A., Godo, L., Torrens, A. and Verdú, V.

**Abstract**

Let K be a variety of (commutative, integral) residuated lattices. The substructural logic usually associated with K is an algebraizable logic that has K as its equivalent algebraic semantics, and is a logic that preserves truth, i.e., 1 is the only truth value preserved by the inferences of the logic. In this paper we introduce another logic associated with K, namely the logic that preserves degrees of truth, in the sense that it preserves lower bounds of truth values in inferences. We study this second logic mainly from the point of view of abstract algebraic logic. We determine its algebraic models and we classify it in the Leibniz and the Frege hierarchies: we show that it is always fully selfextensional, that for most varieties K it is non-protoalgebraic, and that it is algebraizable if and only K is a variety of generalized Heyting algebras, in which case it coincides with the logic that preserves truth. We also characterize the new logic in three ways: by a Hilbert style axiomatic system, by a Gentzen style sequent calculus, and by a set of conditions on its closure operator. Concerning the relation between the two logics, we prove that the truth preserving logic is the extension of the one that preserves degrees of truth with either the rule of Modus Ponens or the rule of Adjunction for the fusion connective.

## 1 Introduction and outline of the paper

This paper concerns (many-valued) logics associated with varieties $\mathsf{K}$ of residuated lattices that, for the sake of simplicity, we assume to be commutative and integral (as usual the maximum of any algebra will be denoted by $1$). These varieties form a wide class of algebraic structures encompassing Heyting algebras, $\mathrm{FL}_{ew}$ algebras, BL algebras, MV algebras, MTL algebras, product algebras, etc. Most of the substructural and fuzzy logics studied in the literature in association with these classes of algebras (see e.g. [40, 41, 44, 45]) are defined by taking 1 as the only truth value to be preserved by inference (in the sense of yielding true consequences from true premises for each interpretation) and they are sometimes called assertional logics [28]. In the present paper they will be referred to as *truth preserving logics* and will be denoted by $\vdash_{\mathsf{K}}$. Truth preserving logics $\vdash_{\mathsf{K}}$ are algebraizable in the sense of Blok and Pigozzi [8], and thus there is a nice correspondence between varieties of residuated lattices and axiomatic extensions of the truth preserving logic of residuated lattices.

The intended role played above by residuated lattices corresponds to truth-value structures, but indeed truth preserving logics do not take full advantage of being many-valued, as they focus on the truth value 1 (*the truth*) and not on other intermediate truth values. A way to circumvent this possible weakness while keeping the truth preserving framework is to introduce truth constants into the language. This methodology goes back to Pavelka [55] where he built a propositional many-valued logic which turned out to be equivalent to the expansion of Łukasiewicz logic by adding a truth constant $\overline{r}$ into the language for each *real* $r \in [0,1]$, together with a number of additional axioms. In this way the expanded language allows one to have formulas of the kind $\overline{r} \to \varphi$ which, when evaluated to 1, express that the truth value of $\varphi$ is greater or equal than $r$. This logic was further developed by Nóvak [51] and Hájek [45]; more recently a similar approach has been applied in [24, 27] to study the expansions with truth constants of other fuzzy logics including Gödel, product, nilpotent minimum logics as well as other continuous t-norm based logics. All these expansions, like the truth preserving logics $\vdash_{\mathsf{K}}$, have been shown to be algebraizable.

In this paper we aim at going beyond the truth preserving framework in order to exploit many-valuedness by focussing on the notion of inference $\models_{\mathsf{K}}^{\leqslant}$ that results from preserving lower bounds of truth values, and hence not only preserving the value 1. In this setting the language remains the same as in truth preserving logics $\vdash_{\mathsf{K}}$, and what changes is the inference relation. This kind of inference corresponds to the so-called *logics preserving degrees of truth* discussed at length in [28, 33, 52], and follows a very general pattern which could be considered for any class of truth structures endowed with an ordering relation. In particular, for a class $\mathsf{K}$ of ordered algebras this pattern would yield the following definition:

$$\begin{aligned} \Gamma \models_{\mathsf{K}}^{\leqslant} \psi \iff & \forall \boldsymbol{A} \in \mathsf{K}\,, \forall v \in \mathrm{Hom}(\boldsymbol{Fm}, \boldsymbol{A})\,, \forall a \in A\,, \\ & \text{if } \ v(\gamma) \geqslant a \text{ for all } \gamma \in \Gamma \ \text{ then } \ v(\psi) \geqslant a. \end{aligned} \tag{1}$$

For each variety $\mathsf{K}$ of residuated lattices, the logics $\models_{\mathsf{K}}^{\leqslant}$ and $\vdash_{\mathsf{K}}$ clearly have the same theorems (i.e., the same consequences from the empty set). We show that the two logics coincide precisely when $\mathsf{K}$ is a variety of generalized Heyting algebras, in other words, when in the algebras of $\mathsf{K}$ the fusion $\star$ coincides with the meet $\wedge$. Notice that this situation corresponds to the cases where the logics might not be considered as properly substructural. Thus, the cases where studying the logics preserving degrees of truth $\models_{\mathsf{K}}^{\leqslant}$ seems to have a genuine interest are the properly substructural ones. In these cases the results in the paper show that the two logics have very different properties, extending the study of Łukasiewicz's logic preserving degrees of truth done in [33].

These differences are apparent in the classification we obtain of the logic $\models_{\mathsf{K}}^{\leqslant}$ with respect to the two main hierarchies in abstract algebraic logic, the Frege and Leibniz hierachies, which deal respectively with replacement properties and with the behaviour of Leibniz operator.

Regarding the Frege hierarchy, while $\models_{\mathsf{K}}^{\leqslant}$ is fully selfextensional for every $\mathsf{K}$, $\vdash_{\mathsf{K}}$ is selfextensional if and only if $\mathsf{K}$ is a variety of generalized Heyting algebras (that is, when $\vdash_{\mathsf{K}}$ coincides with $\models_{\mathsf{K}}^{\leqslant}$; hence it is fully selfextensional).

Moreover, $\vdash_{\mathsf{K}}$ is (fully) Fregean if and only if $\mathsf{K}$ is a variety of generalized Heyting algebras, and these are also the only cases where $\models_{\mathsf{K}}^{\leqslant}$ is (fully) Fregean. As regards the Leibniz hierarchy, it is well known that $\vdash_{\mathsf{K}}$ is always algebraizable (and so protoalgebraic). We prove that $\models_{\mathsf{K}}^{\leqslant}$ is algebraizable if and only if it coincides with $\vdash_{\mathsf{K}}$. Thus, the logics preserving degrees of truth do not give new examples of algebraizabillity. Moreover, we show that $\models_{\mathsf{K}}^{\leqslant}$ is protoalgebraic if and only if there is some $n \in \omega$ such that all algebras in $\mathsf{K}$ satisfy the equation

$$x \wedge \big((x \to y)^n \star (y \to x)^n\big) \preccurlyeq y\,, \qquad (\mathrm{Prot}_n)$$

and that these are precisely the cases where $\models_{\mathsf{K}}^{\leqslant}$ is finitely equivalential. As a consequence, we obtain some new protoalgebraic logics, but we also see that most of the best-known varieties of residuated lattices, such as those mentioned at the beginning, yield natural examples of non-protoalgebraic logics.

The paper also deals with Hilbert style axiomatizability. It is well known that each $\vdash_{\mathsf{K}}$ has an axiomatization with Modus Ponens as the only rule. We show that $\models_{\mathsf{K}}^{\leqslant}$ is axiomatized by using the same axioms as $\vdash_{\mathsf{K}}$, a restricted form of Modus Ponens and the rule of Adjunction for $\wedge$. It is remarkable that by either replacing $\wedge$ by $\star$ in the rule of Adjunction or by using full Modus Ponens, we obtain an axiomatization of $\vdash_{\mathsf{K}}$. As regards the problem of finding axiomatizations exclusively given by Tarski style conditions (in the sense of the term as coined by Wójcicki in [63, p. 107]), it is known that there is no Tarski style axiomatization for $\vdash_{\mathsf{K}}$ [10]. However, as a good approximation we show (Corollary 5.10.3) that $\models_{\mathsf{K}}^{\leqslant}$ can be axiomatized in a way which is very close to a Tarski style axiomatization. To this end, we first give a description of the algebras in the class $\mathsf{K}$ in terms of abstract properties of the closure operator of filter generation (Theorem 5.1), and then use this to define a Gentzen style calculus for $\models_{\mathsf{K}}^{\leqslant}$; this calculus does not seem to be particularly interesting from a proof-theoretic point of view, but the study of its algebraic models yields a new way of relating the class $\mathsf{K}$ to the corresponding logic $\models_{\mathsf{K}}^{\leqslant}$. All this is done explicitly for $\mathsf{K} = \mathsf{RL}$, and then we indicate how to extend the results to other varieties $\mathsf{K}$ of residuated lattices.

The logics we study are related to substructural logics. These logics were initially characterized in connection with proof-theoretic issues [22], but in recent years they have been extensively studied from an algebraic point of view [40]. Indeed, in his survey paper [53], Ono shows how the root of substructurality is closely related to the interaction between fusion and implication in residuated structures, and reaches the conclusion that "substructural logics are the logics of residuated lattices"; of course he has in mind only the truth-preserving logics $\vdash_{\mathsf{K}}$. On the other hand, the results in the present paper establish deep relations between each of the logics $\models_{\mathsf{K}}^{\leqslant}$ and the corresponding class $\mathsf{K}$ of residuated lattices. In a certain sense one could say that they are also logics "of" residuated lattices, although in a different way, and one grcould go even further and conclude that they might also be counted among substructural logics, even if they are not those traditionally considered in the substructural logic literature. This issue is discusssed more in depth in [31, 32]. It is also interesting to look at the alternative definition of substructural logics suggested in [62]; we refer to it in relation with some of our results at the end of Section 4.

The study of consequence relations associated with a class of ordered algebras in the way our $\models_{\mathsf{K}}^{\leqslant}$ are has been rare, but not without precedents, although those specific to the field of many-valued and substructural logics have been very fragmentary or simply exploratory. In his well-known [60, 61], Scott suggests "to replace many values by many valuations", but a detailed analysis of his proposal (see [32]) reveals that he is actually proposing an inference like our $\models^{\leqslant}_{\mathsf{K}}$, in the particular cases corresponding to Łukasiewicz's finitely-valued logics. The same cases where studied in Gil's dissertation [42], with the tools of many-sided sequent calculi. In the case of Gödel's logic (where $\mathsf{K}$ is the variety of linear Heyting algebras), Hájek proves in [45, Theorem 4.2.18] a completeness result of $\vdash_{\mathsf{K}}$ with respect to a semantics (called by him "partial truth") that uses the rational points in the real unit interval in a way very similar to (1). A more explicit formulation is in Baaz, Preining and Zach's [5] (see also [6]), again only for Gödel's logic and only for the semantics based on the real unit interval, but this time for a first-order logic: In this paper the two semantical entailments are explicitly defined and in their Proposition 13 they are proved to coincide (and the remark is made that this does not happen for Łukasiewicz logics). As to the general theory, the term "logics preserving degrees of truth" is coined in Nowak's [52], and in the particular case where the ordered algebras are actually semilattices Jansana's [48] introduces and studies the so-called "semilattice-based logics" with the tools of abstract algebraic logic; some of his results have been applied in the present paper.

The paper is organized as follows. Section 2 presents the logic $\models^{\leqslant}_{\mathsf{K}}$, some of its basic properties, relationships with $\vdash_{\mathsf{K}}$, and a Hilbert style axiomatization. In Section 3 we determine its algebraic models, in particular we show that its algebraic counterpart according to the two general criteria of abstract algebraic logic is indeed $\mathsf{K}$, and that $\models^{\leqslant}_{\mathsf{K}}$ has Leibniz filters. In Section 4 we determine the position of the logic $\models^{\leqslant}_{\mathsf{K}}$ in the Leibniz and Frege hierarchies. In Section 5 we characterize the variety $\mathsf{K}$ in terms of Tarski style conditions (about the closure operator of lattice filter generation) and in terms of Gentzen style rules. Appendix A contains the detailed proofs of several technical lemmas used in the proof of Theorem 5.1. Finally Appendix B exhibits the residuated lattices used as examples and counterexamples in Section 4.

**Acknowledgements.** The collaborative work between the authors of this paper was made possible by several research grants: MTM2008-01139 and TIN2007-68005-C04 of the Spanish Ministry of Education and Science, including FEDER funds of the European Union, and 2005SGR-00083 and 2005SGR-00093 of the Catalan Government. The authors also wish to thank some anonymous referees for helpful comments.

# 2 Definitions, basic properties and axiomatization

We fix the following propositional (algebraic) language $\mathcal{L} = \langle \wedge, \vee, \star, \rightarrow, 1, 0 \rangle$ with connectives of arity $(2, 2, 2, 2, 0, 0)$. We consider only algebras of this similarity type, thus expressions such as "for all algebras" or "an arbitrary algebra" should be understood as restricted to this type. Operations interpreted in a specific algebra will be denoted with the same symbol. All the algebras we will consider have a maximum element, which will be denoted by $1$, and will be the interpretation of the constant 1. By contrast, in the general case we do not presuppose any property of the interpretation of the constant 0. Some more bits of notation: $\boldsymbol{Fm}$ denotes the algebra of formulas (of this similarity type), and $Fm$ its universe. If $\boldsymbol{A}, \boldsymbol{B}$ are two algebras, $\mathrm{Hom}(\boldsymbol{A}, \boldsymbol{B})$ denotes the set of all homomorphisms from $\boldsymbol{A}$ to $\boldsymbol{B}$. If $\boldsymbol{A}$ is an algebra, $\mathrm{Co}\boldsymbol{A}$ denotes the set of all congruences of $\boldsymbol{A}$, and for each class $\mathsf{K}$ of algebras, $\mathrm{Co}_{\mathsf{K}}\boldsymbol{A}$ denotes the subset of $\mathrm{Co}\boldsymbol{A}$ consisting of the congruences $\theta$ of $\boldsymbol{A}$ such that $\boldsymbol{A}/\theta \in \mathsf{K}$.

The largest class of algebras we consider is the class $\mathsf{RL}$ of commutative, integral residuated lattices, which we will call simply ***residuated lattices*** for brevity. Thus, $\boldsymbol{A} = \langle A, \wedge, \vee, \star, \to, 1, 0\rangle \in \mathsf{RL}$ if and only if the reduct $\langle A, \wedge, \vee, 1\rangle$ is a lattice with maximum 1 (its order is denoted by $\leqslant$), the reduct $\langle A, \star, 1\rangle$ is a commutative monoid, and the *fusion* operation $\star$ (sometimes also called the *intensional conjunction* or *strong conjunction*) is residuated, $\to$ being its residual; that is, for all $a, b, c \in A$,

$$a \star b \leqslant c \iff b \leqslant a \to c. \tag{2}$$

It is well-known that the class $\mathsf{RL}$ is a variety; the following forms an equational base for $\mathsf{RL}$:

$$x \wedge y \approx y \wedge x \qquad x \vee y \approx y \vee x \tag{3}$$

$$x \wedge (y \wedge z) \approx (x \wedge y) \wedge z \qquad x \vee (y \vee z) \approx (x \vee y) \vee z \tag{4}$$

$$x \wedge (x \vee y) \approx x \qquad x \vee (x \wedge y) \approx x \tag{5}$$

$$x \wedge 1 \approx x \tag{6}$$

$$x \star (y \star z) \approx (x \star y) \star z \tag{7}$$

$$x \star y \approx y \star x \tag{8}$$

$$x \star 1 \approx x \tag{9}$$

$$x \star (y \vee z) \approx (x \star y) \vee (x \star z) \tag{10}$$

$$x \star (x \to y) \preccurlyeq y \tag{11}$$

$$y \preccurlyeq x \to \big((x \star y) \vee z\big) \tag{12}$$

We use the symbol $\approx$ to represent formal equations, and the symbol $\preccurlyeq$ to represent formal ordering relations; this symbol is interpreted by the order relation of the algebra. Since we deal with lattices, as equations (3)–(5) witness, the expression $\varphi \preccurlyeq \psi$ is actually equivalent to the equation $\varphi \wedge \psi \approx \varphi$.

Notice that residuated lattices may be presented either with or without the constant 0 in the type; in this we adhere to the practice of [40] (see the notes on pages 95 and 97). Having a constant 0 in the type, even if nothing is postulated about it, makes it formally possible to consider several well-known varieties widely studied in substructural and fuzzy logic literature as subvarieties of our $\mathsf{RL}$, which in this setting coincides with the variety of $\mathrm{FL}_{ei}$-algebras of [40]. Nevertheless, the most studied subvarieties of $\mathsf{RL}$ are in fact subvarieties of $\mathsf{FL_{ew}}$, which is the subvariety of $\mathsf{RL}$ obtained when we specify that 0 is the minimum element. Some interesting subvarieties, which we will use later on, are: $\mathsf{MTL}$, the subvariety of $\mathsf{FL_{ew}}$ obtained by adding the prelinearity condition $(x \to y) \vee (y \to x) \approx 1$; $\mathsf{BL}$, the subvariety of $\mathsf{MTL}$ obtained by adding the divisibility condition $x \star (x \to y) \approx x \wedge y$; $\mathsf{MV}$, the subvariety of $\mathsf{BL}$ obtained by making the negation involutive ($\neg\neg x \approx x$, where $\neg x := x \to 0$); $\mathsf{\Pi}$, the variety of product algebras, which is the subvariety of $\mathsf{BL}$ obtained by adding the cancellative property $\neg\neg x \to \big((x \to x \star y) \to (y \star \neg\neg y)\big) \approx 1$; and the variety $\mathsf{G}$ of Gödel or Dummett algebras, obtained from $\mathsf{BL}$ by adding idempotency of fusion $x \star x \approx x$. The main sources of information concerning these varieties and related ones, and the logics associated with them, are [12, 25, 40, 44, 45, 46, 53].

Following modern algebraic logic literature (see for instance [17, 38, 63]), in this paper we identify a ***logic*** $L$ with its consequence relation, which can be denoted as $\vdash_L$ or similar symbols[1]. Since this is a relation $\vdash_L \subseteq P(Fm) \times Fm$, we can use the common relational notation and write $\Gamma \vdash_L \varphi$ instead of $\langle \Gamma, \varphi \rangle \in \vdash_L$; this is to be interpreted as "$\varphi$ follows from $\Gamma$ in the logic $L$". The *interderivability relation* associated with $L$ is the *binary* relation between formulas $\dashv\vdash_L$ defined as follows:

$$\varphi \dashv\vdash_L \psi \iff \varphi \vdash_L \psi \text{ and } \psi \vdash_L \varphi.$$

This relation is always an equivalence relation, but need not be a congruence of the algebra of formulas. This depends on the replacement properties satisfied by $L$. We will deal with this issue in Section 4.

Traditionally, the substructural logic associated in the literature with each subvariety $\mathsf{K}$ of $\mathsf{RL}$ is ***the logic*** $\vdash_{\mathsf{K}}$ ***that preserves truth with respect to the class*** $\mathsf{K}$ (where truth is represented by the constant 1). This logic has $\mathsf{K}$ as its algebraic counterpart; more precisely, $\vdash_{\mathsf{K}}$ is a (finitely and regularly) algebraizable logic having $\mathsf{K}$ as its equivalent algebraic semantics in the sense of [8, 17], with defining equation $x \approx 1$ and equivalence formula

$$x \leftrightarrow y := (x \to y) \star (y \to x).$$

From this it follows that $\vdash_{\mathsf{K}}$ is the (finitary) logic determined by the two clauses

$$\begin{aligned} \varphi_0, \dots, \varphi_{n-1} \vdash_{\mathsf{K}} \psi &\iff \mathsf{K} \models \varphi_0 \approx 1 \;\&\; \dots \;\&\; \varphi_{n-1} \approx 1 \longrightarrow \psi \approx 1, \\ \emptyset \vdash_{\mathsf{K}} \psi &\iff \mathsf{K} \models \psi \approx 1, \end{aligned} \tag{13}$$

where the symbol $\models$, when written without any sub- or superscript, stands for first-order (or quasi-equational) validity, and $\&$ and $\longrightarrow$ are the symbols for conjunction and implication between first-order formulas; the first clause amounts to saying that

$$\begin{aligned} \varphi_0, \dots, \varphi_{n-1} \vdash_{\mathsf{K}} \psi \iff &\forall \boldsymbol{A} \in \mathsf{K}, \forall v \in \mathrm{Hom}(\boldsymbol{Fm}, \boldsymbol{A}), \\ &\text{if } v(\varphi_i) = 1 \text{ for all } i < n, \text{ then } v(\psi) = 1. \end{aligned} \tag{14}$$

A key property of $\vdash_{\mathsf{K}}$ that will be used in the paper is the ***local deduction theorem*** (LDT), that is, the property that for all $\Gamma \cup \{\varphi, \psi\} \subseteq Fm$,

$$\Gamma, \varphi \vdash_{\mathsf{K}} \psi \iff \exists n < \omega \text{ such that } \Gamma \vdash_{\mathsf{K}} \varphi^n \to \psi, \tag{15}$$

where we have used the exponential to abbreviate iterated fusion; that is, we define $x^0 := 1$ and $x^{n+1} := x^n \star x$ for all $n < \omega$.

Following the discussion in the Introduction, we associate with each subvariety $\mathsf{K}$ of $\mathsf{RL}$ another logic. These logics will be our main object of study.

[1]Our usage of the symbols $\vdash$ and $\models$, with different sub- and super-scripts, to denote the consequence relation of a sentential logic is only dictated by practical, even visual reasons. In particular, we do *not* assign, as is often done in the literature, a "syntactical" meaning to $\vdash$ and a "semantical" meaning to $\models$. The only exception is when $\models$ is used as the validity relation of first-order logic, as in the right-hand sides of each of the equivalences in (13).

**Definition 2.1.** *Let* $\mathsf{K}$ *be a subvariety of* $\mathsf{RL}$. *The logic* $\models^{\leqslant}_{\mathsf{K}}$, *which we call* ***the logic that preserves degrees of truth with respect to*** $\mathsf{K}$, *is defined as follows, for all* $\Gamma \cup \{\psi\} \subseteq Fm$:

**1.** *For a finite, non-empty* $\Gamma = \{\varphi_0, \dots, \varphi_{n-1}\}$,

$$\varphi_0, \dots, \varphi_{n-1} \models^{\leqslant}_{\mathsf{K}} \psi \iff \forall \boldsymbol{A} \in \mathsf{K}, \forall v \in \mathrm{Hom}(\boldsymbol{Fm}, \boldsymbol{A}), \forall a \in A, \text{ if } v(\varphi_i) \geqslant a \text{ for all } i < n \text{ then } v(\psi) \geqslant a.$$

**2.** $\emptyset \models^{\leqslant}_{\mathsf{K}} \psi$ *when for all* $\boldsymbol{A} \in \mathsf{K}$, *for all* $v \in \mathrm{Hom}(\boldsymbol{Fm}, \boldsymbol{A})$, $v(\psi) = 1$.

**3.** *For an infinite* $\Gamma \subseteq Fm$, $\Gamma \models^{\leqslant}_{\mathsf{K}} \psi$ *when there exist* $\varphi_0, \dots, \varphi_{n-1} \in \Gamma$ *such that* $\varphi_0, \dots, \varphi_{n-1} \models^{\leqslant}_{\mathsf{K}} \psi$

It follows from clauses 3 and 1 that $\models^{\leqslant}_{\mathsf{K}}$ really *preserves degrees of truth with respect to* $\mathsf{K}$, that is, that it satisfies, for all $\Gamma \cup \{\psi\} \subseteq Fm$,

$$\Gamma \models^{\leqslant}_{\mathsf{K}} \psi \implies \forall \boldsymbol{A} \in \mathsf{K}, \forall v \in \mathrm{Hom}(\boldsymbol{Fm}, \boldsymbol{A}), \forall a \in A, \text{ if } v(\gamma) \geqslant a \text{ for all } \gamma \in \Gamma \text{ then } v(\psi) \geqslant a.$$

It is also easy, reasoning by cases, to see that Definition 2.1 actually yields a (finitary) logic, that is, a finitary consequence relation on the set of formulas.

**Lemma 2.2.** *For all* $\{\varphi_0, \dots, \varphi_{n-1}, \psi\} \subseteq Fm$ *with* $n > 0$,

$$\begin{aligned} \varphi_0, \dots, \varphi_{n-1} \models^{\leqslant}_{\mathsf{K}} \psi &\iff \varphi_0 \wedge \cdots \wedge \varphi_{n-1} \models^{\leqslant}_{\mathsf{K}} \psi \\ &\iff \mathsf{K} \models \varphi_0 \wedge \cdots \wedge \varphi_{n-1} \preccurlyeq \psi \\ &\iff \models^{\leqslant}_{\mathsf{K}} \varphi_0 \wedge \cdots \wedge \varphi_{n-1} \to \psi. \end{aligned}$$

*Proof.* The first equivalence follows by applying to clause 1 of Definition 2.1 the fact that, in a lattice, $v(\varphi_i) \geqslant a$ for all $i < n$ if and only if $v(\varphi_0) \wedge \cdots \wedge v(\varphi_{n-1}) \geqslant a$, if and only if $v(\varphi_0 \wedge \cdots \wedge \varphi_{n-1}) \geqslant a$. The second equivalence is just a rewriting, taking again clause 1 into account. Finally, the third one follows from the fact that in a integral residuated lattice, $a \leqslant b$ is equivalent to $a \to b = 1$ ∎

The first equivalence in Lemma 2.2 is often paraphrased by saying that the logic $\models^{\leqslant}_{\mathsf{K}}$ is ***conjunctive***, as it expresses that the connective $\wedge$ has the behaviour of a conjunction in the classical, extensional sense. Proposition 2.9 will give an algebraic version of this: the filters of this logic on an algebra $\boldsymbol{A} \in \mathsf{K}$ will coincide with the lattice filters of $\boldsymbol{A}$. Note that the equivalence between the first and the last expressions in Lemma 2.2 can be read as a kind of restricted deduction theorem (called graded deduction theorem in [37, Definition 5.3]), where the premisses have to be transferred collectively, in a single step, to the right-hand side of the inference relation.

Now, the second equivalence in Lemma 2.2 together with clause 2 in Definition 2.1 make it clear that *the logic* $\models^{\leqslant}_{\mathsf{K}}$ *depends only on the equations satisfied by the variety* $\mathsf{K}$. Thus, we can actually define it starting from any class $\mathsf{K}$ of residuated lattices, even from a single algebra, using the expression in Lemma 2.2 instead of clause 1 of Definition 2.1, since the logic we obtain is the same as that obtained for the variety generated by $\mathsf{K}$; this will be specially useful when working with examples. For the general theory, however, it is best to work with varieties, hence unless we say otherwise, in the general results about an unspecified class $\mathsf{K}$, we implicitly assume that $\mathsf{K}$ is an arbitrary subvariety of $\mathsf{RL}$.

As a particular case of Lemma 2.2 we have that, for all $\varphi, \psi \in Fm$,

$$\varphi \models_{\mathsf{K}}^{\leqslant} \psi \iff \mathsf{K} \models \varphi \preccurlyeq \psi\,. \tag{16}$$

A direct consequence of this is the following interesting fact about the relation $=\!\!|\!\models_{\mathsf{K}}^{\leqslant}$ of interderivability with respect the consequence $\models_{\mathsf{K}}^{\leqslant}$:

**Corollary 2.3.** *For all* $\varphi, \psi \in Fm$, $\varphi =\!\!|\!\models_{\mathsf{K}}^{\leqslant} \psi \iff \mathsf{K} \models \varphi \approx \psi$. ∎

The logic $\models_{\mathsf{K}}^{\leqslant}$ can be characterized in terms of matrices and generalized matrices. By a ***matrix*** we understand a pair $\langle \boldsymbol{A}, F\rangle$ where $\boldsymbol{A}$ is an algebra and $F$ is a subset of $A$, and by a ***generalized matrix*** (g-matrix for short) we understand a pair $\langle \boldsymbol{A}, \mathcal{C}\rangle$ where $\boldsymbol{A}$ is an algebra and $\mathcal{C}$ is a closed-set system of subsets of $A$ (that is, a family of subsets of $A$ closed under arbitrary intersections and containing $A$). The logic $\vdash_L$ defined by the matrix $\langle \boldsymbol{A}, F\rangle$ is obtained by putting, for all $\Gamma \cup \{\varphi\} \subseteq Fm$,

$$\begin{aligned} \Gamma \vdash_L \varphi \iff & \forall v \in \mathrm{Hom}(\boldsymbol{Fm}, \boldsymbol{A})\,, \\ & \text{if } v(\gamma) \in F\ \forall \gamma \in \Gamma\,,\ \text{then } v(\varphi) \in F. \end{aligned} \tag{17}$$

A matrix $\langle \boldsymbol{A}, F\rangle$ is a ***model*** of a logic $L$ when the implication $\Rightarrow$ in (17) is satisfied; then $F$ is called a ***filter*** of $L$ or $L$***-filter***. For an algebra $\boldsymbol{A}$ we denote by $\mathcal{F}i_{\vdash_L}(\boldsymbol{A})$ the family of all $L$-filters on $\boldsymbol{A}$.

The logic $\vdash_L$ defined by the g-matrix $\langle \boldsymbol{A}, \mathcal{C}\rangle$ is obtained by putting

$$\begin{aligned} \Gamma \vdash_L \varphi \iff & \forall v \in \mathrm{Hom}(\boldsymbol{Fm}, \boldsymbol{A})\,,\ \forall F \in \mathcal{C}, \\ & \text{if } v(\gamma) \in F\ \forall \gamma \in \Gamma\,,\ \text{then } v(\varphi) \in F. \end{aligned} \tag{18}$$

A g-matrix $\langle \boldsymbol{A}, \mathcal{C}\rangle$ is a ***generalized model*** (g-model for short) of $L$ when the implication $\Rightarrow$ in (18) is satisfied.

The logic determined by a class of matrices is defined as the intersection of the logics defined by all the matrices in the family; and similarly for g-matrices. Then:

**Proposition 2.4.** $\models_{\mathsf{K}}^{\leqslant}$ *coincides with the logic defined by the class of matrices*

$$\big\{\langle \boldsymbol{A}, F\rangle : \boldsymbol{A} \in \mathsf{K} \text{ and } F \text{ is a lattice filter of } \boldsymbol{A}\big\}, \tag{19}$$

*and also with the logic defined by the class of g-matrices*

$$\big\{\big\langle \boldsymbol{A}, \mathcal{F}i_{\wedge}(\boldsymbol{A})\big\rangle : \boldsymbol{A} \in \mathsf{K}\big\}, \tag{20}$$

*where* $\mathcal{F}i_{\wedge}(\boldsymbol{A})$ *denotes the family of all the lattice filters of* $\boldsymbol{A}$.

*Proof.* The class of matrices (19) is the union of the bundles of matrices corresponding to the g-matrices in the class (20), therefore the two classes define the same logic. Since $\mathsf{K}$ is a variety and the notion of a lattice filter is elementary definable, the class of matrices (19) is closed under ultraproducts, therefore such logic is finitary. Thus, we need only check that it coincides with $\models^{\leqslant}_{\mathsf{K}}$ for finite sets of assumptions. This follows from the fact that in a lattice, if $Fi_{\wedge}(a_1, \dots, a_n)$ denotes the smallest lattice filter that contains $a_1, \dots, a_n$, then we have that $b \in Fi_{\wedge}(a_1, \dots, a_n)$ if and only if $a_1 \wedge \dots \wedge a_n \leqslant b$, and that $Fi_{\wedge}(\emptyset) = \{1\}$ since our lattices have a maximum $1$. ∎

We now study the relations between the logics $\models^{\leqslant}_{\mathsf{K}}$ and $\vdash_{\mathsf{K}}$. Later on we will need the following easy result, where, according to the previously introduced notation, we denote by $\mathcal{F}i_{\models^{\leqslant}_{\mathsf{K}}}(\boldsymbol{A})$ the set of $\models^{\leqslant}_{\mathsf{K}}$-filters on $\boldsymbol{A}$:

**Lemma 2.5.** *For any $\boldsymbol{A} \in \mathsf{K}$ and every $F \subseteq A$, if $F \in \mathcal{F}i_{\models^{\leqslant}_{\mathsf{K}}}(\boldsymbol{A})$, then $1 \in F$.*

*Proof.* Clause 2 in Definition 2.1 implies that the constant $1$ is a theorem of the logic $\models^{\leqslant}_{\mathsf{K}}$, hence its interpretation in every algebra belongs to every filter of the logic. ∎

Comparing clause 2 of Definition 2.1 with (13) it is straightforward to see that:

**Lemma 2.6.** *For all $\varphi \in Fm$, $\emptyset \models^{\leqslant}_{\mathsf{K}} \varphi$ if and only if $\emptyset \vdash_{\mathsf{K}} \varphi$.* ∎

Thus the two logics *have the same theorems.* This will be a key fact in finding a Hilbert style presentation of $\models^{\leqslant}_{\mathsf{K}}$. Now we investigate the relation between the consequence relations $\models^{\leqslant}_{\mathsf{K}}$ and $\vdash_{\mathsf{K}}$ for non-empty sets of assumptions. First, from Lemmas 2.2 and 2.6 it follows:

**Proposition 2.7.** *For all $\{\varphi_0, \dots, \varphi_{n-1}\} \cup \{\psi\} \subseteq Fm$,*

$$\{\varphi_0, \dots, \varphi_{n-1}\} \models^{\leqslant}_{\mathsf{K}} \psi \iff \emptyset \vdash_{\mathsf{K}} \varphi_0 \wedge \dots \wedge \varphi_{n-1} \to \psi.$$ ∎

Thus, $\models^{\leqslant}_{\mathsf{K}}$ is determined by the theorems of $\vdash_{\mathsf{K}}$. Moreover, this equivalence provides a polynomial reduction of the finitary relation of consequence of $\models^{\leqslant}_{\mathsf{K}}$ to the theoremhood of $\vdash_{\mathsf{K}}$; thus, the complexity of finite theories of $\models^{\leqslant}_{\mathsf{K}}$ is the same as that of the theorems of $\vdash_{\mathsf{K}}$; for certain fuzzy logics this complexity is known, see [3]. It is also interesting to notice that the way in which this reduction is effected coincides with the proposal put forward by Wójcicki [63, Section 2.10] for a general way of associating a logic (a consequence operation) with a set of logical formulas satisfying certain minimal conditions (what he calls an "entailment system"). In addition we have:

**Proposition 2.8.** *The logic $\vdash_{\mathsf{K}}$ is an extension of the logic $\models^{\leqslant}_{\mathsf{K}}$.*

*Proof.* By finitarity and Lemma 2.6, we only need to prove that for any finite non-empty $\{\varphi_0, \dots, \varphi_{n-1}\} \subseteq Fm$ and any $\psi \in Fm$, if $\varphi_0, \dots, \varphi_{n-1} \models^{\leqslant}_{\mathsf{K}} \psi$ then $\varphi_0, \dots, \varphi_{n-1} \vdash_{\mathsf{K}} \psi$. This follows directly from clause 1 of Definition 2.1 and (13), given that $1$ is the maximum of all algebras in $\mathsf{K}$. ∎

Now we are going to find the relations between the filters of the two logics. First of all, by Proposition 2.8 it readily follows that for an arbitrary algebra $\boldsymbol{A}$,

$$\mathcal{F}i_{\vdash_{\mathsf{K}}}(\boldsymbol{A}) \subseteq \mathcal{F}i_{\models^{\leqslant}_{\mathsf{K}}}(\boldsymbol{A}). \tag{21}$$

If $\boldsymbol{A} \in \mathsf{K}$ this relation can be made more precise, as we can characterize both families of logical filters as certain kinds of algebraic filters. Recall that a subset $F$ of a residuated lattice is an ***implicative filter*** when $1 \in F$ and it is closed under Modus Ponens (MP), that is, if $a \in F$ and $a \to b \in F$ then $b \in F$. All implicative filters are lattice filters. It is well-known that the $\vdash_{\mathsf{K}}$-filters in algebras of $\mathsf{K}$, which are residuated lattices, coincide with the implicative filters.

**Proposition 2.9.** *If $\boldsymbol{A} \in \mathsf{K}$ then the filters of the logic $\models_{\mathsf{K}}^{\leqslant}$ on $\boldsymbol{A}$ coincide with the lattice filters of $\boldsymbol{A}$.*

*Proof.* By Proposition 2.4 all lattice filters are filters of the logic. The converse also holds, that is, all filters of the logic will be lattice filters: The property (16) implies that all the filters of the logic will be order filters (increasing subsets); and Lemma 2.2 implies that they will be closed under conjunction. ∎

**Lemma 2.10.** *Let $\boldsymbol{A} \in \mathsf{K}$, and let $F \subseteq A$ be a non-empty order filter (or increasing subset). Then the following conditions are equivalent:*

**1.** *$F$ is an implicative filter.*

**2.** *$F$ is closed under $\star$, that is, if $a, b \in F$ then $a \star b \in F$, for all $a, b \in A$.*

**3.** *$F$ is a lattice filter such that for all $a \in A$, if $a \in F$ then $a^2 = a \star a \in F$.*

*Proof.* $(1 \Rightarrow 2)$ because $a \to (b \to a \star b)) = 1$. $(2 \Rightarrow 1)$ since from $a \star b \leqslant a \wedge b$ and $a \star (a \to b) \leqslant b$, we have that $F$ is a lattice filter closed under MP. $(2 \Rightarrow 3)$ since from $a \star b \leqslant a \wedge b$ we obtain that $F$ is a lattice filter and obviously $a \in F$ implies $a^2 \in F$. $(3 \Rightarrow 2)$ because $(a \wedge b)^2 \leqslant a \star b$. ∎

Since both logics are defined from the class $\mathsf{K}$, they are characterized by their filters on the algebras in $\mathsf{K}$. Therefore:

**Corollary 2.11.** *$\vdash_{\mathsf{K}}$ is the extension of $\models_{\mathsf{K}}^{\leqslant}$ obtained by adding to it any one of the following rules:*

(MP) $\varphi, \varphi \to \psi \ \vdash \ \psi$.

(Adj$\star$) $\varphi, \psi \ \vdash \ \varphi \star \psi$.

(square-closing) $\varphi \ \vdash \ \varphi^2$. ∎

In this statement, by "extension of a logic by adding a rule" we mean the weakest logic containing the given one and satisfying the stated rule. This expression is more clear in the context of Hilbert-style presentations of the logics. We finish the section by finding one for $\models_{\mathsf{K}}^{\leqslant}$.

For each variety $\mathsf{K}$ of residuated lattices, we denote by *Thm*$(\mathsf{K})$ the set of theorems of $\vdash_{\mathsf{K}}$, which by Lemma 2.6 is also the set of theorems of $\models_{\mathsf{K}}^{\leqslant}$. This set of formulas is semantically determined, because it corresponds to the equations of the form $\varphi \approx 1$ true in $\mathsf{K}$, and does not depend on any particular axiomatization of $\vdash_{\mathsf{K}}$. Actually, it is known that $\vdash_{\mathsf{K}}$ can be axiomatized by taking all the formulas in *Thm*$(\mathsf{K})$ as axioms and the rule of Modus Ponens.

Formally, we take an *inference rule* to be any set of pairs $\langle \Gamma, \varphi \rangle$ where $\Gamma$ is a finite set of formulas (if the set can be described by a single scheme then it is commonly written as $\Gamma \vdash \varphi$ or in fraction form as $\dfrac{\Gamma}{\varphi}$). To find an axiomatization of the logic $\models_{\mathsf{K}}^{\leqslant}$ we consider the following two rules of inference:

(Adj-$\wedge$) $\{\langle\{\varphi,\psi\},\varphi\wedge\psi\rangle : \varphi,\psi\in Fm\}$.

(MP-r) $\{\langle\{\varphi,\varphi\to\psi\},\psi\rangle : \varphi,\psi\in Fm \text{ and } \varphi\to\psi\in \mathit{Thm}(\mathsf{K})\}$.

Notice that the rule (MP-r), a restricted form of Modus Ponens, can be applied only when the premise of the form $\varphi\to\psi$ belongs to $\mathit{Thm}(\mathsf{K})$.

**Theorem 2.12** (Completeness). *The logic* $\models^{\leqslant}_{\mathsf{K}}$ *is axiomatized by taking all the formulas in* $\mathit{Thm}(\mathsf{K})$ *as axioms, and the rules* (Adj-$\wedge$) *and* (MP-r) *as rules of inference.*

*Proof.* Let us first see that each axiom and rule is satisfied by $\models^{\leqslant}_{\mathsf{K}}$:

– If $\varphi\in\mathit{Thm}(\mathsf{K})$, then by Lemma 2.6 $\emptyset\models^{\leqslant}_{\mathsf{K}}\varphi$.

– The rule (Adj-$\wedge$): Lemma 2.2 in particular implies that $\varphi,\psi\models^{\leqslant}_{\mathsf{K}}\varphi\wedge\psi$.

– The rule (MP-r): If $\varphi\to\psi\in\mathit{Thm}(\mathsf{K})$ and $v\in\mathrm{Hom}(\boldsymbol{Fm},\boldsymbol{A})$ for some $\boldsymbol{A}\in\mathsf{K}$, then $v(\varphi\to\psi)=v(\varphi)\to v(\psi)=1$, so $v(\varphi)\leqslant v(\psi)$ and then $v(\varphi)\wedge v(\varphi\to\psi)=v(\varphi)\leqslant v(\psi)$. By Lemma 2.2 again, this means that $\{\varphi,\varphi\to\psi\}\models^{\leqslant}_{\mathsf{K}}\psi$.

Now, proceeding by induction on the length of the proof it is easy to see that, if there is a proof of $\psi$ from a set of assumption $\varphi_0,\dots,\varphi_{n-1}$ in the axiomatic system, then $\varphi_0,\dots,\varphi_{n-1}\models^{\leqslant}_{\mathsf{K}}\psi$.

Conversely, assume now that $\varphi_0,\dots,\varphi_{n-1}\models^{\leqslant}_{\mathsf{K}}\psi$; by Proposition 2.7 this means that $\emptyset\vdash_{\mathsf{K}}(\varphi_0\wedge\cdots\wedge\varphi_{n-1})\to\psi$, so $(\varphi_0\wedge\cdots\wedge\varphi_{n-1})\to\psi\in\mathit{Thm}(\mathsf{K})$ is an axiom of the system. Therefore, $n-1$ applications of the rule (Adj-$\wedge$) followed by one application of the rule (MP-r) yield a proof of $\psi$ from $\{\varphi_0,\dots,\varphi_{n-1}\}$ in the axiomatic system. This ends the completeness proof. ∎

If, as it happens in most of the examples studied in the literature, one has an axiomatic presentation for $\vdash_{\mathsf{K}}$ that uses a set of axioms $AX(\mathsf{K})$ and only the ordinary rule of Modus Ponens, then the set $AX(\mathsf{K})$ can replace the set $\mathit{Thm}(\mathsf{K})$ in the presentation given in Theorem 2.12. This is so because each application of Modus Ponens in a proof of a theorem of $\vdash_{\mathsf{K}}$ will in fact be an application of (MP-r), therefore the same proof will be a proof in $\models^{\leqslant}_{\mathsf{K}}$.

The algebraizability of $\vdash_{\mathsf{K}}$ implies that, from any equational presentation of the variety $\mathsf{K}$, an axiomatic presentation of the logic $\vdash_{\mathsf{K}}$ can be obtained, see [9, Theorem 8.0.9]. This presentation has Modus Ponens among its rules, but has other inference rules besides the axioms (basically, the rules expressing congruence of the equivalence formulas). To find an axiomatic presentation of $\vdash_{\mathsf{K}}$ using only Modus Ponens as rule, which would yield an axiomatic presentation of $\models^{\leqslant}_{\mathsf{K}}$, one can use the local deduction theorem (15) and this will turn all these rules into theorems, which can then be taken as axioms.

When the set of theorems of $\vdash_{\mathsf{K}}$ is decidable, both the set of axioms and the set of rules of our axiomatization of $\models^{\leqslant}_{\mathsf{K}}$ will be decidable, and in this case the given axiomatic presentation could be called "in the Hilbert style" in full sense of the term. Clearly, this will happen when the equational theory of $\mathsf{K}$ is decidable. Notice that, even when the set of theorems of $\vdash_{\mathsf{K}}$ is not decidable but recursively enumerable, using Craig's technique [16] it is not hard to replace our proposed axioms and rules for $\models^{\leqslant}_{\mathsf{K}}$ by other ones that are decidable.

After Theorem 2.12, Corollary 2.11 can be read, more or less informally, as saying that $\models^{\leqslant}_{\mathsf{K}}$ is obtained from $\vdash_{\mathsf{K}}$ by weakening the rules of Modus Ponens and Adjunction, namely, restricting the rule of Modus Ponens to (MP-r) and replacing Adjunction for $\star$ by Adjunction for $\wedge$.

# 3 Algebras and models

As we already pointed out, the variety $\mathsf{K}$ is the equivalent algebraic semantics of the logic $\vdash_{\mathsf{K}}$. This means that $\mathsf{K}$ is *the algebraic counterpart* of $\vdash_{\mathsf{K}}$ according to a well established general paradigm of the algebraization of a logic, namely the theory of *algebraizable logics* initiated by Blok and Pigozzi in [8] (see also [17]). As to the relation between $\mathsf{K}$ and the logic $\models_{\mathsf{K}}^{\leqslant}$, it goes beyond the completeness result contained implicitly in its definition, but its description needs a more general framework because, as we will show in Section 4, in general the logic $\models_{\mathsf{K}}^{\leqslant}$ is not algebraizable.

We are going to consider two general, abstract theories of the algebraization of a logic. Both define the algebraic counterpart of a logic starting from the notion of a reduced model, which is defined in terms of the Leibniz congruence. If $\langle \boldsymbol{A}, F\rangle$ is an arbitrary matrix, its ***Leibniz congruence*** $\boldsymbol{\Omega}_{\boldsymbol{A}}F$ is the largest of all congruences of $\boldsymbol{A}$ that are ***compatible*** with $F$ in the sense that they do not relate elements inside $F$ with elements outside $F$; that is,

$$\boldsymbol{\Omega}_{\boldsymbol{A}}F = \max\big\{\theta \in \mathrm{Co}\boldsymbol{A} : \text{ if } \langle a,b\rangle \in \theta \text{ and } a \in F \text{ then } b \in F\big\}.$$

It can be shown that such a congruence always exists and can be characterized syntactically. A matrix is ***reduced*** when $\boldsymbol{\Omega}_{\boldsymbol{A}}F$ is the identity; this means that the identity relation is the only congruence of the algebra that is compatible with $F$. Notice that these notions are purely algebraic ones, and are not influenced by the fact that the matrix is or is not a model of some logic.

The first abstract theory we consider is the general theory of matrices [63], which forms the oldest part of abstract algebraic logic. This theory takes as *algebraic conterpart of a logic* $L$ the class $\mathsf{Alg}^*(L)$ of algebra reducts of the reduced models of $L$:

$$\mathsf{Alg}^*(L) = \big\{\boldsymbol{A} : \exists F \subseteq A\,, \langle \boldsymbol{A}, F\rangle \text{ is a model of } L \text{ and } \boldsymbol{\Omega}_{\boldsymbol{A}}F = Id\big\}.$$

If $L$ is algebraizable, then $\mathsf{Alg}^*(L)$ coincides with its equivalent algebraic semantics. This is also the class considered by Rasiowa for the so-called *implicative logics* treated in [59].

A more general theory is introduced in [35], it is proposed to take generalized models as the tool to associate a class of algebras with a logic. If $\langle \boldsymbol{A}, \mathcal{C}\rangle$ is a g-matrix, then its ***Tarski congruence*** is the largest congruence of $\boldsymbol{A}$ that is compatible with all the members of $\mathcal{C}$; it can be seen equal to

$$\widetilde{\boldsymbol{\Omega}}_{\boldsymbol{A}}\mathcal{C} = \bigcap_{F\in\mathcal{C}} \boldsymbol{\Omega}_{\boldsymbol{A}}F.$$

A g-matrix is ***reduced*** when its Tarski congruence is the identity, and then the class of $L$***-algebras*** $\mathsf{Alg}(L)$ is defined as the class of algebra reducts of the reduced g-models of $L$, that is,

$$\mathsf{Alg}(L) = \big\{\boldsymbol{A} : \exists \mathcal{C} \subseteq P(A)\,, \langle \boldsymbol{A}, \mathcal{C}\rangle \text{ is a g-model of } L \text{ and } \widetilde{\boldsymbol{\Omega}}_{\boldsymbol{A}}\mathcal{C} = Id\big\}.$$

It is easy to see that $\mathsf{Alg}(L)$ is the class of subdirect products of algebras in $\mathsf{Alg}^*(L)$, so in general $\mathsf{Alg}^*(L) \subseteq \mathsf{Alg}(L)$. The empirical study of examples and the theoretical work developed in [35] and subsequent papers (see [29, 38]) supports the adoption of the class $\mathsf{Alg}(L)$ as a better definition of the algebraic counterpart of a logic $L$ in the most general case.

The two theories agree in a large class of logics, first considered by Blok and Pigozzi in [7]: the class of *protoalgebraic logics*, which are those where the behaviour of matrix semantics and of the Leibniz operator is good enough to obtain a smooth transfer between logical and algebraic properties. In [35, Proposition 3.2] it is shown that for a protoalgebraic logic $L$, $\mathsf{Alg}^*(L) = \mathsf{Alg}(L)$, but in general these two classes can be different. Since, as we show in Section 4, most of the logics $\models_{\mathsf{K}}^{\leqslant}$ are not protoalgebraic, we should identify the class $\mathsf{Alg}(\models_{\mathsf{K}}^{\leqslant})$ in order to determine the algebraic counterpart of $\models_{\mathsf{K}}^{\leqslant}$.

To this end, notice that by Lemma 2.2, $\models_{\mathsf{K}}^{\leqslant}$ is ***semilattice based*** with respect to $\mathsf{K}$ in the sense of [48, p. 76]. This enables us to draw a number of consequences. To begin with, for a semilattice based logic $L$, Proposition 3.1 of [48] contains a practical characterization of the class $\mathsf{Alg}(L)$: If $L$ is an arbitrary semilattice based logic with respect to $\mathsf{K}$ then $\mathsf{Alg}(L)$ is the so-called ***intrinsic variety*** of $L$, that is, the class of algebras (variety) defined by the equations $\varphi \approx \psi$ such that $\varphi$ and $\psi$ are interderivable with respect to $L$, that is, such that $\varphi \dashv\vdash_L \psi$. In our case, where $L$ is $\models_{\mathsf{K}}^{\leqslant}$, $\mathsf{Alg}(\models_{\mathsf{K}}^{\leqslant})$ is determined by $=\!\!|\!\models_{\mathsf{K}}^{\leqslant}$, and then Corollary 2.3 yields:

**Proposition 3.1.** $\mathsf{Alg}(\models_{\mathsf{K}}^{\leqslant}) = \mathsf{K}$. ∎

In the theory of [35], among the g-models of a logic a prominent role is played by the so-called ***full g-models***, i.e., the g-models $\langle \boldsymbol{A}, \mathcal{C} \rangle$ of a logic $L$ that are of the form $\mathcal{C} = \{h^{-1}[F] : F \in \mathcal{F}i_{\vdash_L}(\boldsymbol{B})\}$ for certain homomorphism $h$ from $\boldsymbol{A}$ onto $\boldsymbol{B}$. In Corollary 2.12 of [35] it is shown that $\mathsf{Alg}(L)$ is also the class of algebra reducts of the reduced full g-models of $L$; in our case these models can be characterized in a very precise way:

**Proposition 3.2.** *Let $\langle \boldsymbol{A}, \mathcal{C} \rangle$ be an arbitrary g-matrix. Then $\langle \boldsymbol{A}, \mathcal{C} \rangle$ is a reduced full g-model of $\models_{\mathsf{K}}^{\leqslant}$ if and only if $\boldsymbol{A} \in \mathsf{K}$ and $\mathcal{C} = \mathcal{F}i_{\wedge}(\boldsymbol{A})$.*

*Proof.* If $\langle \boldsymbol{A}, \mathcal{C} \rangle$ is a reduced full g-model of $\models_{\mathsf{K}}^{\leqslant}$, in particular it is a reduced g-model, and by the general definition of $\mathsf{Alg}(L)$ mentioned above, $\boldsymbol{A} \in \mathsf{Alg}(\models_{\mathsf{K}}^{\leqslant})$, that is, $\boldsymbol{A} \in \mathsf{K}$ by Proposition 3.1. Now we can apply Theorem 2.30 of [35], which implies that on each algebra there is a unique reduced full g-model of a logic. Since by definition the g-matrix $\langle \boldsymbol{A}, \mathcal{F}i_{\models_{\mathsf{K}}^{\leqslant}}(\boldsymbol{A}) \rangle$ is always a full g-model, together with Proposition 2.9 this implies that $\mathcal{C} = \mathcal{F}i_{\wedge}(\boldsymbol{A})$. This also proves the converse. ∎

In turn, this yields a characterization of the class $\mathsf{K}$:

**Corollary 3.3.** *Let $\boldsymbol{A}$ be any algebra. Then $\boldsymbol{A} \in \mathsf{K}$ if and only if there is a family $\mathcal{C}$ of filters of the logic $\models_{\mathsf{K}}^{\leqslant}$ such that the g-matrix $\langle \boldsymbol{A}, \mathcal{C} \rangle$ is reduced.* ∎

Proposition 3.2 shows that we can always take the family of all lattice filters of the algebra $\boldsymbol{A}$ as the family $\mathcal{C}$.

The general study of semilattice based logics in [48] does not determine $\mathsf{Alg}^*(L)$ in general; but here we can determine it by an ad-hoc argument:

**Proposition 3.4.** $\mathsf{Alg}^*(\models_{\mathsf{K}}^{\leqslant}) = \mathsf{K}$.

*Proof.* Since $\mathsf{Alg}^*(\models_{\mathsf{K}}^{\leqslant}) \subseteq \mathsf{Alg}(\models_{\mathsf{K}}^{\leqslant})$ is valid in general, by Proposition 3.1 we have that $\mathsf{Alg}^*(\models_{\mathsf{K}}^{\leqslant}) \subseteq \mathsf{K}$. To prove the converse inclusion, take any $\boldsymbol{A} \in \mathsf{K}$. Then by the algebraizability of $\vdash_{\mathsf{K}}$ with defining equation $x \approx 1$, the matrix $\langle \boldsymbol{A}, \{1\} \rangle$ is a reduced model of $\vdash_{\mathsf{K}}$, and by Proposition 2.8 every model of $\vdash_{\mathsf{K}}$ is also a model of $\models_{\mathsf{K}}^{\leqslant}$. Therefore $\langle \boldsymbol{A}, \{1\} \rangle$ is a reduced model of $\models_{\mathsf{K}}^{\leqslant}$, and by the definition of $\mathsf{Alg}^*(L)$, this implies that $\boldsymbol{A} \in \mathsf{Alg}^*(\models_{\mathsf{K}}^{\leqslant})$. ∎

Note that Propositions 3.1 and 3.4 show that $\mathsf{K}$ is the algebraic counterpart of $\models_{\mathsf{K}}^{\leqslant}$ in a very strong sense, as the two general theories coincide in this case.

The particular case of this result for $\mathsf{K} = \mathsf{MV}$, the variety of MV-algebras, is stated in Corollary 3.5 of [33]. Unfortunately, the proof given there is flawed, but we have now proved it more in general.

We end the section with a technical result which will be needed in the next section and which will allow us to obtain a characterization of the reduced models of $\models_{\mathsf{K}}^{\leqslant}$, something that can not be obtained from the general theory:

**Proposition 3.5.** *For all* $\boldsymbol{A}$ *and all* $F \in \mathcal{F}i_{\models_{\mathsf{K}}^{\leqslant}}(\boldsymbol{A})$*, there exists the set*

$$F^+ = \min\{G \in \mathcal{F}i_{\models_{\mathsf{K}}^{\leqslant}}(\boldsymbol{A}) : \boldsymbol{\Omega}_{\boldsymbol{A}} G = \boldsymbol{\Omega}_{\boldsymbol{A}} F\}. \tag{22}$$

*This set* $F^+$ *satisfies the following properties:*

1. $F^+$ *is a* $\vdash_{\mathsf{K}}$*-filter.*
2. $F^+ \subseteq F$.
3. $F^+$ *is the only* $\vdash_{\mathsf{K}}$*-filter having* $\boldsymbol{\Omega}_{\boldsymbol{A}} F$ *as its Leibniz congruence.*
4. $F^+ = 1/\boldsymbol{\Omega}_{\boldsymbol{A}} F$.

*Proof.* From the general definitions it follows that if $F \in \mathcal{F}i_{\models_{\mathsf{K}}^{\leqslant}}(\boldsymbol{A})$ then $\boldsymbol{\Omega}_{\boldsymbol{A}} F \in \mathrm{Co}_{\mathsf{Alg}^*(\models_{\mathsf{K}}^{\leqslant})} \boldsymbol{A}$, that is, $\boldsymbol{\Omega}_{\boldsymbol{A}} F$ is a congruence of $\boldsymbol{A}$ such that the quotient $\boldsymbol{A}/\boldsymbol{\Omega}_{\boldsymbol{A}} F$ is in $\mathsf{Alg}^*(\models_{\mathsf{K}}^{\leqslant})$. By Proposition 3.4, $\mathsf{Alg}^*(\models_{\mathsf{K}}^{\leqslant}) = \mathsf{K}$, thus $\boldsymbol{\Omega}_{\boldsymbol{A}} F \in \mathrm{Co}_{\mathsf{K}} \boldsymbol{A}$. On the other hand, since the stronger logic $\vdash_{\mathsf{K}}$ is algebraizable with equivalent algebraic semantics $\mathsf{K}$, the mapping $\boldsymbol{\Omega}_{\boldsymbol{A}}$ is an isomorphism between the lattices $\mathcal{F}i_{\vdash_{\mathsf{K}}}(\boldsymbol{A})$ and $\mathrm{Co}_{\mathsf{K}} \boldsymbol{A}$ (see [38, Theorem 3.13] for instance); thus, there exists a unique $F^+ \in \mathcal{F}i_{\vdash_{\mathsf{K}}}(\boldsymbol{A})$ such that $\boldsymbol{\Omega}_{\boldsymbol{A}} F^+ = \boldsymbol{\Omega}_{\boldsymbol{A}} F$. Since the equation defining the algebraizability of $\vdash_{\mathsf{K}}$ is $x \approx 1$, $F^+ = 1/\boldsymbol{\Omega}_{\boldsymbol{A}} F^+ = 1/\boldsymbol{\Omega}_{\boldsymbol{A}} F$ (which shows 4). We prove that $F^+ \subseteq F$: if $a \in F^+$ then $\langle a, 1 \rangle \in \boldsymbol{\Omega}_{\boldsymbol{A}} F$, but by Lemma 2.5, we have that $1 \in F$ and by the compatibility of $\boldsymbol{\Omega}_{\boldsymbol{A}} F$ with $F$, $a \in F$. This shows 2, and the points 1 and 3 are true by construction. By (21), $F^+ \in \mathcal{F}i_{\models_{\mathsf{K}}^{\leqslant}}(\boldsymbol{A})$, therefore $F^+ \in \{G \in \mathcal{F}i_{\models_{\mathsf{K}}^{\leqslant}} : \boldsymbol{\Omega}_{\boldsymbol{A}} G = \boldsymbol{\Omega}_{\boldsymbol{A}} F\}$ and we have to prove that it is its minimum: take any $G \in \{G \in \mathcal{F}i_{\models_{\mathsf{K}}^{\leqslant}}(\boldsymbol{A}) : \boldsymbol{\Omega}_{\boldsymbol{A}} G = \boldsymbol{\Omega}_{\boldsymbol{A}} F\}$ and consider the set $G^+$ constructed in the same way as $F^+$: We already know that $G^+ \subseteq G$, that $G^+ \in \mathcal{F}i_{\vdash_{\mathsf{K}}}(\boldsymbol{A})$, and that $\boldsymbol{\Omega}_{\boldsymbol{A}} G^+ = \boldsymbol{\Omega}_{\boldsymbol{A}} G$. Therefore, $\boldsymbol{\Omega}_{\boldsymbol{A}} G^+ = \boldsymbol{\Omega}_{\boldsymbol{A}} G = \boldsymbol{\Omega}_{\boldsymbol{A}} F = \boldsymbol{\Omega}_{\boldsymbol{A}} F^+$. By the mentioned isomorphism, this implies $F^+ = G^+$, so $F^+ \subseteq G$. ∎

The filters with the property (22) were introduced and studied in depth in [36], where they are called ***Leibniz filters***. There it is shown that they exist for every protoalgebraic logic, but their existence for non-protoalgebraic logics is unknown in general. This gives an independent interest to Proposition 3.5. Moreover, it is easy to use it to obtain the following characterization of the reduced models of $\models_{\mathsf{K}}^{\leqslant}$.

**Corollary 3.6.** *A matrix $\langle \boldsymbol{A}, F\rangle$ is a reduced model of $\models_{\mathsf{K}}^{\leqslant}$ if and only if $\boldsymbol{A} \in \mathsf{K}$ and $F$ is a lattice filter of $\boldsymbol{A}$ such that $F^+ = \{1\}$.*

*Proof.* If $\langle \boldsymbol{A}, F\rangle$ is a reduced model of $\models_{\mathsf{K}}^{\leqslant}$, by Proposition 3.4 $\boldsymbol{A} \in \mathsf{K}$ and by Proposition 2.9 $F$ is a lattice filter of $\boldsymbol{A}$. By algebraizability, $\{1\}$ is the only filter of $\vdash_{\mathsf{K}}$ on $\boldsymbol{A}$ having the identity as its Leibniz congruence. Hence $F^+ = \{1\}$. The converse is also obvious given that $\boldsymbol{\Omega}_{\boldsymbol{A}} F = \boldsymbol{\Omega}_{\boldsymbol{A}} F^+$. ∎

There are several ways to rephrase the condition that $F^+ = \{1\}$. For instance, by part 3 of Proposition 3.5, $F^+ = \{1\}$ is equivalent to saying that $\{1\}$ is the only filter $G$ of $\vdash_{\mathsf{K}}$ whose Leibniz congruence $\boldsymbol{\Omega}_{\boldsymbol{A}} G$ is compatible with $F$. Or, by part 4 of Proposition 3.5, this is to say that if $\langle a, 1\rangle \in \boldsymbol{\Omega}_{\boldsymbol{A}} F$ then $a = 1$. However, we do not have a working characterization of $\boldsymbol{\Omega}_{\boldsymbol{A}} F$ when $F \in \mathcal{F}i_{\models_{\mathsf{K}}^{\leqslant}}(\boldsymbol{A})$ for an arbitrary $\boldsymbol{A}$; thus the characterization of the reduced matrices in Corollary 3.6 is of little use in general. However, in the next section we will find better ones in case the logic is protoalgebraic, and indeed some of them will be shown to be equivalent to the property of being so.

# 4 Classification in the hierarchies of abstract algebraic logic

In this section we investigate the location of the logic $\models_{\mathsf{K}}^{\leqslant}$ inside the two general hierarchies of logics considered in abstract algebraic logic: the Leibniz hierarchy and the Frege hierarchy; see [17, 30, 38] for more details on the hierarchies and on abstract algebraic logic in general. The Leibniz hierarchy is organized around properties of the Leibniz operator on the filters of the logics; we deal with it later. Let us begin with the Frege hierarchy, which is organized, so to speak, around the kind of replacement properties of the logic. The basic, largest class in this hierarchy is that of ***selfextensional logics***, which are the logics $L$ whose interderivability relation $\dashv\vdash_L$ is a congruence of the algebra of formulas. This means that $L$ satisfies the following *weak form of the replacement property*: For any $\alpha, \beta, \varphi(x, \vec{y}) \in Fm$,

$$\text{if} \quad \alpha \dashv\vdash_L \beta \quad \text{then} \quad \varphi(\alpha, \vec{y}) \dashv\vdash_L \varphi(\beta, \vec{y}). \tag{23}$$

We see that all our logics belong to this class:

**Proposition 4.1.** *For each variety $\mathsf{K}$ of residuated lattices, the consequence relation $\models_{\mathsf{K}}^{\leqslant}$ is selfextensional.*

*Proof.* We have seen in Corollary 2.3 that $\varphi \mathrel{=\!\models}_{\mathsf{K}}^{\leqslant} \psi$ holds if and only if $\mathsf{K} \models \varphi \approx \psi$. This last relation is always a congruence relation. ∎

The next class in the hierarchy is that of ***fully selfextensional*** logics, which are the logics whose full g-models inherit the property of replacement (23); for this it is enough to postulate that for each algebra $\boldsymbol{A} \in \mathsf{Alg}(L)$, the interderivability relation of the "basic" full g-model on $\boldsymbol{A}$, which is $\langle \boldsymbol{A}, \mathcal{F}i_{\vdash_L} \boldsymbol{A} \rangle$, is a congruence of $\boldsymbol{A}$. Then:

**Corollary 4.2.** *For each variety* $\mathsf{K}$ *of residuated lattices, the consequence relation* $\models^{\leqslant}_{\mathsf{K}}$ *is fully selfextensional.*

*Proof.* By Proposition 3.1, for each $\mathsf{K}$ we have that $\mathsf{Alg}(\models^{\leqslant}_{\mathsf{K}}) = \mathsf{K}$ and by Proposition 2.9, if $\boldsymbol{A} \in \mathsf{K}$ the filters of $\models^{\leqslant}_{\mathsf{K}}$ on $\boldsymbol{A}$ are the lattice filters of $\boldsymbol{A}$. It is well-kown that in a lattice, if two elements generate the same filter then they are equal. Thus the interderivability relation associated with the g-matrix $\big\langle \boldsymbol{A}, \mathcal{F}i_{\models^{\leqslant}_{\mathsf{K}}}(\boldsymbol{A})\big\rangle$ is the identity, which is trivially a congruence. Since by Proposition 3.2 the g-matrices of this form are the reduced full g-models of $\models^{\leqslant}_{\mathsf{K}}$, this establishes that this logic is fully selfextensional. ∎

Fully selfextensional logics form a smaller class than that of selfextensional logics in the Frege hierarchy, and they have better properties. The other two classes in the Frege hierarchy are those of Fregean logics and of fully Fregean ones. They are defined similarly to the preceding ones, but requiring the postulated replacement properties to hold relatively to every theory of the logic (and of the full g-models, respectively). Thus, a logic $L$ is ***Fregean*** when for each theory $T$ of $L$, the interderivability relation of $L$ modulo $T$ is a congruence of $\boldsymbol{Fm}$. This corresponds to the following *strong form of the replacement property*: For any theory $T$ of $L$, any $\alpha, \beta, \varphi(x, \vec{y}) \in Fm$,

$$\text{if} \quad T, \alpha \vdash_L \beta \;\text{ and }\; T, \beta \vdash_L \alpha \quad \text{then} \quad T, \varphi(\alpha, \vec{y}) \vdash_L \varphi(\beta, \vec{y}). \tag{24}$$

Fregean logics were introduced in [35, 56], and have been extensively studied in the context of protoalgebraic logics in [19, 20], see also [17]. ***Fully Fregean logics*** are those that inherit the Fregean character for all their full g-models, and form the smallest and best behaved class in the hierarchy. Next we see that among our logics only those corresponding to varieties of generalized Heyting algebras belong to this class, and that they coincide with those that are Fregean. ***Generalized Heyting algebras*** can be informally described as Heyting algebras possibily without minimum; relatively to the variety $\mathsf{RL}$ they are defined by the equation

$$x \star y \approx x \wedge y. \tag{25}$$

These algebras are called *relatively pseudo-complemented lattices* in [59], and *Brouwerian algebras*[2] in [40] and other works. Had we developed our theory in the restricted case of subvarieties of $\mathsf{FL_{ew}}$, then "Heyting algebras" should replace "generalized Heyting algebras" everywhere in the paper.

**Proposition 4.3.** *Let* $\mathsf{K}$ *be a variety of residuated lattices. Then the following conditions are equivalent:*

**1.** *The logic* $\models^{\leqslant}_{\mathsf{K}}$ *is fully Fregean.*

[2]Notice that the name "Brouwerian algebras" has also been used in the literature to denote the duals of Heyting algebras.

**2.** *The logic* $\models^{\leqslant}_{\mathsf{K}}$ *is Fregean.*

**3.** $\mathsf{K}$ *is a variety of generalized Heyting algebras.*

*Proof.* ($1 \Rightarrow 2$) is trivial.
($2 \Rightarrow 3$) Let us note that each variable $x$ is interderivable with $1$ modulo the theory generated by $x$, that is, $\{x, x\} \models^{\leqslant}_{\mathsf{K}} 1$ and $\{x, 1\} \models^{\leqslant}_{\mathsf{K}} x$. From the assumption that $\models^{\leqslant}_{\mathsf{K}}$ is Fregean we infer that $x \star x$ and $1 \star 1$ are interderivable modulo the same theory. In particular, $\{x, 1 \star 1\} \models^{\leqslant}_{\mathsf{K}} x \star x$, that is, $x \models^{\leqslant}_{\mathsf{K}} x \star x$. Since the converse inference always holds, we conclude that $x \approx x^2$ holds in $\mathsf{K}$, which, as is well known, implies that all the members of $\mathsf{K}$ are generalized Heyting algebras.
($3 \Rightarrow 1$) If $\boldsymbol{A} \in \mathsf{K}$ then by Proposition 2.9 $\mathcal{F}i_{\models^{\leqslant}_{\mathsf{K}}}(\boldsymbol{A}) = \mathcal{F}i_{\wedge}(\boldsymbol{A})$, and it is well known that in a generalized Heyting algebra, for each lattice filter $F$, the relation $a \equiv b \iff Fi_{\wedge}(F, a) = Fi_{\wedge}(F, b)$ is a congruence. This tells us that $\models^{\leqslant}_{\mathsf{K}}$ is fully Fregean. ∎

Summarizing, our class of logics is divided into two groups as far as the Frege hierarchy is concerned: Those generated by a class of generalized Heyting algebras are not properly substructural because they satisfy (25), and are fully Fregean; the rest are fully selfextensional but not Fregean. We are going to find the same division later on in this section.

Now we start the study of the classification of our logics in the Leibniz hierarchy. This hierarchy is much richer and more complicated than the Fregean one, and most of its classes can be defined or characterized in several ways, although the majority of them concern several properties of ***the Leibniz operator***; this is the mapping

$$F \longmapsto \boldsymbol{\Omega}_{\boldsymbol{A}} F$$

restricted to the set of all $L$-filters on a fixed algebra $\boldsymbol{A}$.

The largest class in the Leibniz hierarchy is that of ***protoalgebraic logics***. They can be defined in a number of equivalent ways. We first consider them as the logics such that the monotonicity condition

$$\text{if } G \subseteq F \text{ then } \boldsymbol{\Omega}_{\boldsymbol{A}} G \subseteq \boldsymbol{\Omega}_{\boldsymbol{A}} F \tag{M}$$

holds for every algebra $\boldsymbol{A}$ and every pair of filters $F, G$ of the logic over $\boldsymbol{A}$. One can also consider the following particular case of condition (M):

$$\text{if } G \subseteq F \text{ and } \boldsymbol{\Omega}_{\boldsymbol{A}} F = Id \text{ then } \boldsymbol{\Omega}_{\boldsymbol{A}} G = Id \tag{M'}$$

It is not difficult to see that condition (M') is strictly weaker than (M) and is not sufficient to guarantee protoalgebraicity; one example is the fragment of classical propositional logic with only $\wedge$, $\vee$ and $1$, see [34], and the reason is that the filters of the reduced models of this logic have exactly one point. Next we prove, among other facts, that in our case this restriction is indeed sufficient. In the proof we use the well known fact that if $F$ is a $\vdash_{\mathsf{K}}$-filter then $\boldsymbol{\Omega}_{\boldsymbol{A}} F = \{\langle a, b\rangle \in A \times A : a \leftrightarrow b \in F\}$.

**Theorem 4.4.** *Let* $\mathsf{K}$ *be a variety of residuated lattices. Then the following conditions are equivalent:*

**1.** *The logic* $\models^{\leqslant}_{\mathsf{K}}$ *is protoalgebraic.*

**2.** *For all* $\boldsymbol{A}$ *and all* $F \in \mathcal{F}i_{\models^{\leqslant}_{\mathsf{K}}}(\boldsymbol{A})$*,* $F^{+} = \max\{G \in \mathcal{F}i_{\vdash_{\mathsf{K}}}(\boldsymbol{A}) : G \subseteq F\}$*.*

**3.** *A matrix* $\langle \boldsymbol{A}, F\rangle$ *is a reduced model of* $\models^{\leqslant}_{\mathsf{K}}$ *if and only if* $\boldsymbol{A} \in \mathsf{K}$ *and* $F$ *is a lattice filter of* $\boldsymbol{A}$ *such that* $\{1\}$ *is the only implicative filter of* $\boldsymbol{A}$ *contained in* $F$*.*

**4.** *For every* $\boldsymbol{A} \in \mathsf{K}$ *and every* $F, G \in \mathcal{F}i_{\models^{\leqslant}_{\mathsf{K}}}(\boldsymbol{A})$*, if* $G \subseteq F$ *and* $\boldsymbol{\varOmega_A} F = Id$ *then* $\boldsymbol{\varOmega_A} G = Id$*.*

*Proof.* ($1 \Rightarrow 2$) If $G$ is a $\vdash_{\mathsf{K}}$-filter such that $G \subseteq F$, then it is also a $\models^{\leqslant}_{\mathsf{K}}$-filter so, if $\models^{\leqslant}_{\mathsf{K}}$ is protoalgebraic, $\boldsymbol{\varOmega_A} G \subseteq \boldsymbol{\varOmega_A} F = \boldsymbol{\varOmega_A} F^{+}$. But $G$ and $F^{+}$ are also $\vdash_{\mathsf{K}}$-filters, thus by the algebraizability of $\vdash_{\mathsf{K}}$, the mapping $\boldsymbol{\varOmega_A}$ is an order isomorphism on the lattice of $\vdash_{\mathsf{K}}$-filters, therefore $G \subseteq F^{+}$.

($2 \Rightarrow 3$) If $\boldsymbol{A} \in \mathsf{K}$ then $\{1\}$ is certainly the smallest implicative filter contained in a given lattice filter. The result then follows from Corollary 3.6, taking Propositions 3.4 and 2.9 into account.

($3 \Rightarrow 4$) If $\langle \boldsymbol{A}, F\rangle$ is a reduced model of $\models^{\leqslant}_{\mathsf{K}}$ then the characterization in 3 implies that $\langle \boldsymbol{A}, G\rangle$ will also be a reduced model whenever $G \subseteq F$.

($4 \Rightarrow 1$) Let $\boldsymbol{A}$ be an arbitrary algebra and $F, G \in \mathcal{F}i_{\models^{\leqslant}_{\mathsf{K}}}(\boldsymbol{A})$ with $G \subseteq F$. We are going to prove that $\boldsymbol{\varOmega_A} G \subseteq \boldsymbol{\varOmega_A} F$. By Proposition 3.5 we can replace $G$ with $G^{+}$, and therefore without loss of generality we can assume that $G \in \mathcal{F}i_{\vdash_{\mathsf{K}}}(\boldsymbol{A})$. Let us denote by $G \vee F^{+}$ the smallest $\vdash_{\mathsf{K}}$-filter including $G \cup F^{+}$. It is clear that $F^{+} \subseteq G \cup F^{+} \subseteq G \vee F^{+}$. Next we prove that $G \vee F^{+} \subseteq F$. Since $G$ and $F^{+}$ are $\vdash_{\mathsf{K}}$-filters and $F$ is an $\models^{\leqslant}_{\mathsf{K}}$-filter, in order to prove that $G \vee F^{+} \subseteq F$ it is enough by Lemma 2.10 to take $a \in G$ and $b \in F^{+}$ and show that $a \star b \in F$. Using that $b \leqslant a \leftrightarrow (a \star b)$ it follows that $a \leftrightarrow (a \star b) \in F^{+}$. Therefore, $\langle a, a \star b\rangle \in \boldsymbol{\varOmega_A} F^{+} = \boldsymbol{\varOmega_A} F$. Finally, from the previous fact together with that $a \in G \subseteq F$ we get that $a \star b \in F$.

Thus, we know that $F^{+} \subseteq G \vee F^{+} \subseteq F$. Now we define $\widetilde{\boldsymbol{A}} = \boldsymbol{A}/\boldsymbol{\varOmega_A} F^{+}$, $\widetilde{G} = (G \vee F^{+})/\boldsymbol{\varOmega_A} F^{+}$ and $\widetilde{F} = F/\boldsymbol{\varOmega_A} F^{+}$. By Proposition 3.4 we know that $\widetilde{\boldsymbol{A}} \in \mathsf{K}$. And since $\boldsymbol{\varOmega_A} F = \boldsymbol{\varOmega_A} F^{+}$ we know that $\boldsymbol{\varOmega}_{\widetilde{\boldsymbol{A}}} \widetilde{F} = Id$. Using the assumption it follows that $\boldsymbol{\varOmega}_{\widetilde{\boldsymbol{A}}} \widetilde{G} = Id$. Thus, $\boldsymbol{\varOmega_A} F^{+} = \boldsymbol{\varOmega_A}(G \vee F^{+})$. Using that $\boldsymbol{\varOmega_A}$ is an isomorphism between $\mathcal{F}i_{\vdash_{\mathsf{K}}}(\boldsymbol{A})$ and $\mathrm{Co}_{\mathsf{K}}\boldsymbol{A}$ (because $\vdash_{\mathsf{K}}$ is algebraizable) we obtain that $F^{+} = G \vee F^{+}$, i.e., $G \subseteq F^{+}$. Using again the isomorphism we get that $\boldsymbol{\varOmega_A} G \subseteq \boldsymbol{\varOmega_A} F^{+}$. Hence, since $\boldsymbol{\varOmega_A} F = \boldsymbol{\varOmega_A} F^{+}$ it follows that $\boldsymbol{\varOmega_A} G \subseteq \boldsymbol{\varOmega_A} F$. ∎

Notice that the statement in part 2 of Theorem 4.4 contains two different facts: it says that for a given filter $F$ of $\models^{\leqslant}_{\mathsf{K}}$ there exists the largest filter of $\vdash_{\mathsf{K}}$ contained in it, and it says that this largest filter coincides with $F^{+}$. It is interesting to observe that in general this maximum need not exist, as shown in Example 1 in Appendix B. But notice that a weaker property always holds: by a typical application of Zorn's Lemma, if $G$ is a filter of $\vdash_{\mathsf{K}}$ and $G \subseteq F$ for a filter $F$ of $\models^{\leqslant}_{\mathsf{K}}$, then $G$ can be extended to a filter of $\vdash_{\mathsf{K}}$ that is maximal among those contained in $F$. Even more, the mentioned maximum may exist in non-protoalgebraic logics without being equal to $F^{+}$. One example where this happens is the class of MTL algebras. The maximum exists in each algebra in this class because if $\boldsymbol{A} \in \mathsf{MTL}$ and $F$ is a lattice filter of $\boldsymbol{A}$, then the set $\{a \in A : a^{n} \in F \text{ for every } n < \omega\}$ is an implicative filter (because $a^{2n} \wedge b^{2n} = (a \wedge b)^{2n} = (a \wedge b)^{n} \star (a \wedge b)^{n} \leqslant a^{n} \star b^{n} = (a \star b)^{n}$) and it is indeed the maximum implicative filter inside $F$. On the other hand, as it will be shown later on, the logic $\models^{\leqslant}_{\mathsf{MTL}}$ is not protoalgebraic (see item 1 after Theorem 4.9). Therefore, Theorem 4.4 implies that there must exist a MTL algebra (e.g., Example 2 in Appendix B) with a filter $F$ of $\models^{\leqslant}_{\mathsf{MTL}}$ where the largest filter of $\vdash_{\mathsf{MTL}}$ inside $F$ exists without being equal to $F^{+}$.

Let $\mathsf{MV}$ denote the class of MV algebras, which are the algebraic counterpart of Łukasiewicz's infinite-valued logic. In [33, Theorem 3.11] it was shown that $\models^{\leqslant}_{\mathsf{MV}}$ is non-protoalgebraic. In Example 2 in Appendix B a much simpler proof is given, by considering Chang's algebra. Another interesting algebra where (M) fails is the product algebra given in Example 3 of Appendix B; this proves that $\models^{\leqslant}_{\Pi}$ is non-protoalgebraic, where $\Pi$ is the class of product algebras.

It is obvious that non-protoalgebraic logics must always yield instances of algebras and filters where (M) fails. However, this seems not the best strategy to prove non-protoalgebraicity. The difficulty is that there are varieties (e. g., $\mathsf{MV}$ and $\Pi$) where the Leibniz operator is monotonic over lattice filters in a generator of the variety (e. g., when the generator is a simple algebra) while there are other algebras in the variety where the monotonicity fails. Hence, providing an equational characterization of protoalgebraicity has the advantage over the previous one that it is enough to check it in the algebras generating the variety.

The proof is based on another characterization of protoalgebraicity of $\models^{\leqslant}_{\mathsf{K}}$: the existence of a set of formulas $\varDelta(x,y)$ in two variables such that

$$\emptyset \models^{\leqslant}_{\mathsf{K}} \varDelta(x,x) \quad \text{and} \quad x\,, \varDelta(x,y) \models^{\leqslant}_{\mathsf{K}} y \tag{P}$$

Since $\models^{\leqslant}_{\mathsf{K}}$ is finitary and conjunctive (2.2), the set $\varDelta(x,y)$ can actually be reduced to just one formula $\delta(x,y)$. We are going to see that we can always take $(x \to y)^n \star (y \to x)^n$, for some $n < \omega$, as the formula $\delta(x,y)$. Notice that in general there is no uniqueness, up to equivalence in $\models^{\leqslant}_{\mathsf{K}}$, of the formula $\delta(x,y)$ satisfying (P); for instance, for classical propositional logic we can consider either the formula $x \to y$ or the formula $x \leftrightarrow y$, which are not equivalent in this logic.

Consider, for any natural number $n$, the variety $\mathsf{Prot_n} \subseteq \mathsf{RL}$ of the residuated lattices that satisfy the following equation:

$$x \wedge \big((x \to y)^n \star (y \to x)^n\big) \preccurlyeq y \tag{Prot$_n$}$$

Note that, by commutativity and associativity of $\star$, this is actually the same as $x \wedge (x \leftrightarrow y)^n \preccurlyeq y$. This variety can be alternatively defined by other equations with logical significance:

**Theorem 4.5.** *Let $n \geqslant 1$. Any of the following equations can replace* (Prot$_n$) *in the definition of the variety* $\mathsf{Prot_n}$ *relatively to the variety* $\mathsf{RL}$*:*

$$x \wedge y^n \approx x \star y^n \tag{SC$_n$}$$

$$x \wedge (x \to y)^n \preccurlyeq y \tag{MP$_n$}$$

$$x \wedge (x \to y)^n \wedge (y \to x)^n \preccurlyeq y \tag{SMP$_n$}$$

*Proof.* Let us denote provisionally[3] by $\mathsf{SC_n}, \mathsf{MP_n}, \mathsf{SMP_n}$ the varieties defined by the respective equations. It is straightforward that $\mathsf{SC_n} \subseteq \mathsf{MP_n} \subseteq \mathsf{SMP_n} \subseteq$

[3] The names have been chosen so as to follow the associations SC: "Strong Contraction", MP: "Modus Ponens" and SMP: "Symmetric Modus Ponens".

$\mathsf{Prot_n}$, so we have only to show that $\mathsf{Prot_n} \subseteq \mathsf{SC_n}$ in order to establish the statement. First we show that if $\boldsymbol{A} \in \mathsf{Prot_n}$ then

$$\boldsymbol{A} \models x^n \approx x^{n+1}. \tag{26}$$

If $a \in A$, by residuation $a \leqslant a^n \to a^{n+1} = a^n \leftrightarrow a^{n+1}$. Hence, $a^n \leqslant (a^n \leftrightarrow a^{n+1})^n$, i.e., $a^n = a^n \wedge (a^n \leftrightarrow a^{n+1})^n$. Using ($\mathrm{Prot}_n$) this implies that $a^n \leqslant a^{n+1}$. Thus, $a^n = a^{n+1}$, which shows (26).

Now, in order to prove that $\mathsf{Prot_n} \subseteq \mathsf{SC_n}$, let $\boldsymbol{A} \in \mathsf{Prot_n}$ and let $a, b \in A$. It is enough to prove that $a \wedge b^n \leqslant a \star b^n$, but by (26) it is enough to check that $a \wedge b^{n \cdot n} \leqslant a \star b^n$. By residuation $b^n \leqslant a \to (a \star b^n) = a \leftrightarrow (a \star b^n)$. Hence, $a \wedge b^{n \cdot n} \leqslant a \wedge \big((a \leftrightarrow (a \star b^n)\big)^n$. Using ($\mathrm{Prot}_n$) this implies that $a \wedge b^{n \cdot n} \leqslant a \star b^n$. This completes the proof that $\mathsf{SC_n} = \mathsf{MP_n} = \mathsf{SMP_n} = \mathsf{Prot_n}$. ∎

The simplest and more workable among the four equivalent equations defining the class $\mathsf{Prot_n}$ seems to be condition ($\mathrm{SC}_n$), as it uses only the operations $\wedge$ and $\star$. Recall that the operation $\star$ is used to defined the exponential notation, so that the other conditions use the three operations $\wedge, \star$ and $\to$.

The significance of the family of classes $\mathsf{Prot_n}$, as well as the justification for its name, lies in the following important result:

**Theorem 4.6.** *Let* $\mathsf{K}$ *be a variety of residuated lattices. Then the following conditions are equivalent:*

**1.** *The logic* $\models_{\mathsf{K}}^{\leqslant}$ *is protoalgebraic.*

**2.** *There is an* $n < \omega$ *such that* $\mathsf{K} \models x \wedge \big((x \to y)^n \star (y \to x)^n\big) \preccurlyeq y$, *that is, such that* $\mathsf{K} \subseteq \mathsf{Prot_n}$.

*Proof.* ($2 \Rightarrow 1$) If we take $\delta(x, y) := (x \to y)^n \star (y \to x)^n$ then the first condition in (P) holds trivially, and the second one is equivalent to the assumption.
($1 \Rightarrow 2$) Consider the algebra of formulas $\boldsymbol{Fm}$ and the set

$$\theta := \big\{ \langle \varphi, \psi \rangle \in Fm \times Fm : (x \to y) \star (y \to x) \vdash_{\mathsf{K}} (\varphi \to \psi) \star (\psi \to \varphi) \big\}.$$

It is well known that $\theta$ is a congruence of $\boldsymbol{Fm}$ (it is actually the Leibniz congruence of the $\vdash_{\mathsf{K}}$-theory generated by $(x \to y) \star (y \to x)$), and obviously $\langle x, y \rangle \in \theta$ (indeed, it is the smallest one with this property). By the Local Deduction Theorem (15) for $\vdash_{\mathsf{K}}$ it is easy to see that $\theta$ is compatible with the theory $T$ of $\models_{\mathsf{K}}^{\leqslant}$ generated by the set $\big\{ (x \to y)^n \star (y \to x)^n : n < \omega \big\}$. Hence, $\theta \subseteq \boldsymbol{\Omega}_{\boldsymbol{Fm}} T$. Let $T'$ be the theory of $\models_{\mathsf{K}}^{\leqslant}$ generated by $T \cup \{x\}$. By the monotonicity given by the protoalgebraicity assumption, $\theta \subseteq \boldsymbol{\Omega}_{\boldsymbol{Fm}} T'$, i.e., $\theta$ is compatible with $T'$. Since $\langle x, y \rangle \in \theta$ and $x \in T'$, we conclude that $y \in T'$, that is,

$$\{x\} \cup \big\{ (x \to y)^n \star (y \to x)^n : n < \omega \big\} \models_{\mathsf{K}}^{\leqslant} y$$

Since $\models_{\mathsf{K}}^{\leqslant}$ is finitary, there is an $n < \omega$ such that

$$\emptyset \models_{\mathsf{K}}^{\leqslant} \big( x \wedge \big( (x \to y)^n \star (y \to x)^n \big) \big) \to y,$$

which is equivalent to the statement in 2. ∎

Conditions (26) also define an interesting family of subvarieties, as do other equations which will be of interest for our research. Let us formally introduce them. Consider the equations

$$x \vee (x^n \to y) \approx 1 \tag{EM$_n$}$$

$$x^n \wedge y^n \approx x^n \star y^n \tag{IMC$_n$}$$

$$x^n \approx x^{n+1} \tag{E$_n$}$$

The varieties of residuated lattices given by these equations will be denoted, respectively, by $\mathsf{EM_n}$, $\mathsf{IMC_n}$ and $\mathsf{E_n}$. The varieties $\mathsf{EM_n}$ and $\mathsf{E_n}$ have been widely considered in the literature under the same names, see [40, pp. 96 and 463], and they correspond, respectively, to the variety generated by the simple $n$-contractive residuated lattices and to the variety of the $n$-contractive residuated lattices. As for $\mathsf{IMC_n}$, the name relates to the property called "Idempotent Meet Contraction" (cf. Lemma 4.7).

It is straightforward to check that for every $\mathsf{K} \in \{\mathsf{EM}, \mathsf{Prot}, \mathsf{IMC}, \mathsf{E}\}$, if $1 \leqslant n \leqslant m$ then $\mathsf{K_n} \subseteq \mathsf{K_m}$. Indeed, using the MV chain with $n+2$ points we get that if $1 \leqslant n < m$ then $\mathsf{K_n} \subsetneq \mathsf{K_m}$. Moreover, the cases $n = 0$ and $n = 1$ yield well-known varieties: For $n = 0$ all the previously defined varieties coincide with the trivial variety, except in the cases of $\mathsf{SC}$ and $\mathsf{IMC}$, where $\mathsf{SC_0} = \mathsf{RL} = \mathsf{IMC_0}$. For $n = 1$ it is not difficult to show that $\mathsf{EM_1}$ is the class of generalized Boolean algebras while $\mathsf{Prot_1} = \mathsf{IMC_1} = \mathsf{E_1}$ is the class of generalized Heyting algebras. In case we add the condition that 0 is the minimum then we obtain the usual classes of Boolean and Heyting algebras, respectively (and $\mathrm{FL}_{ew}$-algebras instead of $\mathsf{RL}$).

It is well known that in every residuated lattice the set of idempotent elements is closed under join, and that in general this is not the case for meet, as Example 1 in Appendix B shows. However, it is a simple exercise to show:

**Lemma 4.7.** *Let $n \geqslant 1$. If $\boldsymbol{A} \in \mathsf{E_n}$, then the following conditions are equivalent:*

**1.** $\boldsymbol{A} \models x^n \wedge y^n \approx x^n \star y^n$.

**2.** $\boldsymbol{A} \models x^n \wedge y^n \approx (x \wedge y)^n$.

**3.** *The set of idempotent elements of $\boldsymbol{A}$ is closed under meet.* ∎

It is not difficult to realize (see Lemma 4.8) that we can replace the equation (IMC$_n$) by the equation (EM$_n$) plus any of the properties in the previous lemma. It is obvious that in all MTL algebras idempotent elements are closed under meet, but notice that the class of residuated lattices where idempotent elements are closed under meet is strictly bigger than $\mathsf{MTL}$, as witnessed for instance by Example 4 in Appendix B.

In the proof of Theorem 4.5 we have shown that $\mathsf{Prot_n} \subseteq \mathsf{E_n}$. However we can do better:

**Lemma 4.8.** *Let $n \geqslant 1$. Then,* $\mathsf{EM_n} \subseteq \mathsf{Prot_n} \subseteq \mathsf{IMC_n} \subseteq \mathsf{E_n}$.

*Proof.* Let us start with $\mathsf{EM_n} \subseteq \mathsf{Prot_n}$. It is well known that $\mathsf{EM_n}$ is the variety generated by simple $n$-contractive algebras (see [40, Chapter 11]). Now suppose that $\boldsymbol{A} \models$ (EM$_n$) and $\boldsymbol{A}$ is a simple $n$-contractive residuated latice. Then, for every pair of elements $a, b \in A$ it holds that either $b = 1$ or $b^n = 0$. Thus, in both cases, $a \wedge b^n = a \star b^n$. Therefore, $\boldsymbol{A} \models (\mathrm{SC}_n)$. After Theorem 4.5, this establishes that $\mathsf{EM_n} \subseteq \mathsf{Prot_n}$. The other two inclusions are straightforward. ∎

In the next result we establish, among several things, that the inclusions are proper. We notice that in the literature there are at least two different ways to consider an ordinal sum of a family of residuated lattices: one way [50, 45] identifies the maximum 1 of one algebra with the minimum of the next algebra (this applies only to bounded algebras and when the family is ordered, usually well-ordered), and the other way [2] identifies the maximum 1 of all the algebras in the family. In this paper we always use the term "ordinal sum" in this last sense, which indeed is the standard approach of ordinal sums in the case of hoops and semihoops.

**Theorem 4.9.** *Let $n \geqslant 2$. Then:*

1. *The equations $(\mathrm{Prot}_n), (\mathrm{IMC}_n)$ and $(\mathrm{E}_n)$, and their equivalents, are preserved under the operation of ordinal sums, while $(\mathrm{EM}_n)$ is not.*
2. $\mathsf{EM_n} \subsetneq \mathsf{Prot_n} \subsetneq \mathsf{IMC_n} \subsetneq \mathsf{E_n}$.
3. $\mathsf{MTL} \cap \mathsf{EM_n} \subsetneq \mathsf{MTL} \cap \mathsf{Prot_n} \subsetneq \mathsf{MTL} \cap \mathsf{IMC_n} = \mathsf{MTL} \cap \mathsf{E_n}$.
4. $\mathsf{MV} \cap \mathsf{EM_n} = \mathsf{BL} \cap \mathsf{EM_n} \subsetneq \mathsf{BL} \cap \mathsf{Prot_n} = \mathsf{BL} \cap \mathsf{E_n}$.
5. *If $\boldsymbol{A}$ is an MTL chain, then $\boldsymbol{A} \in \mathsf{Prot_n}$ if and only if $\boldsymbol{A}$ is an ordinal sum of simple $n$-contractive MTL chains.*

*Proof.* 1. The negative part follows from the fact that the ordinal sum of the Łukasiewicz chain of $n+1$-values with itself does not satisfy $(\mathrm{EM}_n)$. For the positive part (which also works for $n \geqslant 0$) we only need to take care about the case where the two elements involved in the equations are in different components, and this can be straightforwardly checked by the reader.

2. The fact that $\mathsf{EM_n} \subsetneq \mathsf{Prot_n}$ is witnessed by the ordinal sum of the Łukasiewicz chain of $n+1$-values with itself; indeed, it shows that $\mathsf{BL} \cap \mathsf{EM_n} \subsetneq \mathsf{BL} \cap \mathsf{Prot_n}$, and hence a fortiori that $\mathsf{BL} \cap \mathsf{Prot_n} \not\subseteq \mathsf{EM_n}$. The strict inclusion $\mathsf{Prot_n} \subsetneq \mathsf{IMC_n}$ is witnessed by Example 5 in Appendix B; actually, it shows that $\mathsf{MTL} \cap \mathsf{Prot_n} \subsetneq \mathsf{MTL} \cap \mathsf{IMC_n}$. Finally, $\mathsf{IMC_n} \subsetneq \mathsf{E_n}$ is witnessed by Example 1 in Appendix B.

3. The same counterexamples that were given in part 2 work here. And the only new inclusion is that $\mathsf{MTL} \cap \mathsf{E_n} \subseteq \mathsf{IMC_n}$. Let $\boldsymbol{A} \in \mathsf{MTL} \cap \mathsf{E_n}$. We want to prove that $\boldsymbol{A} \models x^n \wedge y^n \approx x^n \star y^n$. Since we are considering MTL algebras we can suppose that $\boldsymbol{A}$ is a chain. Let $a, b \in A$. Then, $a^n \wedge b^n = (a \wedge b)^n = (a \wedge b)^n \star (a \wedge b)^n \leqslant a^n \star b^n \leqslant a^n \wedge b^n$. Therefore, all these inequalities are equalities.

4. The first equality follows from [40, Theorem 11.19], [14, Theorem 1.7] and [15, Proposition 3.1]. For the inequality, the same counterexamples that were used in part 2 work here again. Finally, we have to show that $\mathsf{BL} \cap \mathsf{E_n} \subseteq \mathsf{Prot_n}$. We remind the reader that $\mathsf{BL} \models x \wedge y \approx y \star (y \to x)$. Take any $\boldsymbol{A} \in \mathsf{BL} \cap \mathsf{E_n}$ and any $a, b \in A$. Then $a \wedge b^n = b^n \star (b^n \to a) = b^n \star b^n \star (b^n \to a) \leqslant b^n \star a \leqslant a \wedge b^n$. Therefore, $a \wedge b^n = a \star b^n$ and so $\boldsymbol{A} \models (\mathrm{SC}_n)$, that is, $\boldsymbol{A} \in \mathsf{Prot_n}$.

5. By 1 it is trivial that if $\boldsymbol{A}$ is an ordinal sum of simple $n$-contractive MTL chains then $\boldsymbol{A} \models (\mathrm{Prot}_n)$. Let us assume now that $\boldsymbol{A} \models (\mathrm{Prot}_n)$. Then we know that $\boldsymbol{A} \models (\mathrm{E}_n)$. Hence, two elements $a, b \in A$ are in the same

Archimedean component iff $a^n = b^n$. It is obvious that every Archimedean component is a simple $n$-contractive MTL chain. We only need to check that $\boldsymbol{A}$ is the ordinal sum of its Archimedean components. In other words, we have to check that if $a, b$ are two elements in different components then $a \star b = a \wedge b$. Let us assume that $a, b$ are in different components and that $a < b$ (remember that we are in a chain). Since they are in different components it holds that $a < b^n$. Thus, $a = a \wedge b^n \leqslant a \wedge \big(a \to (a \star b)\big)^n \leqslant (a \star b)^n \leqslant a \star b \leqslant a \wedge b \leqslant a$. Therefore, all these inequalities are equalities. ∎

We note that MTL chains coincide with $\mathrm{FL}_{ew}$ chains, and in this setting the restriction to chains in part 5 of Theorem 4.9 is unavoidable; this necessity is witnessed for instance by Example 4 in Appendix B.

As a consequence of Theorem 4.9 we have the following facts:

1. If $\mathsf{K}$ is a variety of residuated lattices such that $\models^{\leqslant}_{\mathsf{K}}$ is protoalgebraic then there exists an $n < \omega$ such that all algebras of $\mathsf{K}$ are $n$-contractive (i.e., $\mathsf{K} \subseteq \mathsf{E_n}$). In particular, it follows that if $\mathsf{K}$ is any of the varieties $\mathsf{RL}, \mathsf{MV}, \mathsf{\Pi}, \mathsf{BL}, \mathsf{MTL}$ or $\mathsf{FL_{ew}}$ then $\models^{\leqslant}_{\mathsf{K}}$ is non-protoalgebraic.
2. If $\mathsf{K}$ is a variety of residuated lattices such that $\models^{\leqslant}_{\mathsf{K}}$ is protoalgebraic then in all algebras of $\mathsf{K}$ the meet of idempotent elements is always an idempotent element.
3. If $\mathsf{K} \subseteq \mathsf{MTL}$, then the logic $\models^{\leqslant}_{\mathsf{K}}$ is protoalgebraic if and only if there exists a natural $n$ such that all chains in $\mathsf{K}$ are ordinal sums of simple $n$-contractive MTL chains. We point out that there are finite MTL chains not defining a protoalgebraic logic, for instance Example 5 in Appendix B.
4. If $\mathsf{K} \subseteq \mathsf{BL}$, then the logic $\models^{\leqslant}_{\mathsf{K}}$ is protoalgebraic if and only if there exists a natural $n$ such that $\mathsf{K} \subseteq \mathsf{E_n}$. In particular any finite BL chain defines a protoalgebraic logic.
5. If $\mathsf{K}$ is the variety generated by a family of continuous t-norms (over the real unit interval), then $\mathsf{K}$ defines a protoalgebraic logic if and only if $\mathsf{K}$ is the variety of Gödel algebras. This is a consequence of the fact that the standard Gödel t-norm is the only continuous t-norm that is $n$-contractive for some $n \geqslant 2$; and indeed it is $n$-contractive for every $n \in \omega$. In fact, the logic preserving the degrees of truth over Gödel algebras is not only protoalgebraic but algebraizable, as we show below.
6. The equation ($\mathrm{Prot}_n$) (and its equivalents) gives a characterization of ordinal sums of simple $n$-contractive MTL chains that is alternative (possibly simpler) than the one stated in [47, Prop. 4.24]. In this paper these ordinal sums were characterized using (both) equations $x^n \approx x^{n+1}$ and $(y^n \to x) \vee \big(x \to (x \star y)\big) \approx 1$.

Now we investigate when the logics preserving degrees of truth are equivalential. By definition $\models^{\leqslant}_{\mathsf{K}}$ is ***equivalential*** when there is a set of formulas $\Delta(x, y)$ in two variables satisfying the condition (P) and the condition

$$\Delta(x, y) \cup \Delta(z, w) \models^{\leqslant}_{\mathsf{K}} \Delta(x \circ z, y \circ w) \text{ where } \circ \in \{\star, \wedge, \vee, \to\}. \tag{E}$$

The formulas in the set $\Delta(x, y)$ are called ***equivalence formulas***. By definition, every equivalential logic is protoalgebraic, but notice that for an equivalen-tial logic, it is not true that any set $\varDelta(x,y)$ satisfying (P), and hence witnessing protoalgebraicity, also satisfies (E). Moreover, in equivalential logics the set $\varDelta(x,y)$ involved in the definition is unique up to interderivability in the logic, something not true for protoalgebraic logics. A logic is ***finitely equivalential*** when it is equivalential and has a finite set of equivalence formulas. In our case, since $\models_{\mathsf{K}}^{\leqslant}$ is finitary and conjunctive, when $\models_{\mathsf{K}}^{\leqslant}$ is finitely equivalential we can assume that $\varDelta(x,y)$ reduces to just one formula $\delta(x,y)$.

**Theorem 4.10.** *Let* $\mathsf{K}$ *be a variety of residuated lattices and let* $n \geqslant 1$. *Then the following conditions are equivalent:*

1. $\mathsf{K} \subseteq \mathsf{Prot_n}$.
2. *For every algebra* $\boldsymbol{A}$ *and every* $F \in \mathcal{F}i_{\models_{\mathsf{K}}^{\leqslant}}(\boldsymbol{A})$, $F^+ = \{a \in A : a^n \in F\}$.
3. *For every algebra* $\boldsymbol{A}$ *and every* $F \in \mathcal{F}i_{\models_{\mathsf{K}}^{\leqslant}}(\boldsymbol{A})$, $\boldsymbol{\varOmega_A} F = \big\{\langle a,b\rangle \in A \times A : (a \leftrightarrow b)^n \in F\big\}$.
4. *The logic* $\models_{\mathsf{K}}^{\leqslant}$ *is finitely equivalential, with* $\varDelta(x,y) := \big\{(x \leftrightarrow y)^n\big\}$ *as set of equivalence formulas.*

*Proof.* $(1 \Rightarrow 2)$ Assume that $\mathsf{K} \subseteq \mathsf{Prot_n}$ and that $F \in \mathcal{F}i_{\models_{\mathsf{K}}^{\leqslant}}(\boldsymbol{A})$, and let $G$ be $\{a \in A : a^n \in F\}$. Since $\mathsf{Prot_n} \subseteq \mathsf{IMC_n}$ it is clear that $G$ is closed under $\star$. Thus, $G$ is a $\vdash_{\mathsf{K}}$-filter of $\boldsymbol{A}$. Indeed, $G$ is the largest $\vdash_{\mathsf{K}}$-filter inside $F$. We obtain that $F^+ = G$ by Theorems 4.4 and 4.6.

$(2 \Rightarrow 3)$ is a consequence of the algebraizability of $\vdash_{\mathsf{K}}$.

$(3 \Rightarrow 4)$ is a general fact in abstract algebraic logic, related to another of the characterizations of equivalential logics: the definability of the Leibniz congruence in the filters of the logic.

$(4 \Rightarrow 1)$ follows from the facts that in general, in an equivalential logic every set $\varDelta(x,y)$ of equivalence formulas also satisfies (P), and that for the logic $\models_{\mathsf{K}}^{\leqslant}$ and the set $\varDelta(x,y) := \big\{(x \leftrightarrow y)^n\big\}$ the second property in P amounts to the condition (Prot$_n$). ∎

If we take the fact that $\mathsf{Prot_n} \subseteq \mathsf{E_n}$ into account, then Theorem 4.10 yields:

**Corollary 4.11.** *Let* $\mathsf{K}$ *be a variety of residuated lattices. Then the following conditions are equivalent:*

1. *The logic* $\models_{\mathsf{K}}^{\leqslant}$ *is protoalgebraic.*
2. *The logic* $\models_{\mathsf{K}}^{\leqslant}$ *is finitely equivalential.*
3. *For every algebra* $\boldsymbol{A}$ *and* $F \in \mathcal{F}i_{\models_{\mathsf{K}}^{\leqslant}}(\boldsymbol{A})$ *it holds that* $F^+ = \{a \in A : a^n \in F \textit{ for every } n < \omega\}$.
4. *For every algebra* $\boldsymbol{A}$ *and* $F \in \mathcal{F}i_{\models_{\mathsf{K}}^{\leqslant}}(\boldsymbol{A})$ *it holds that* $\boldsymbol{\varOmega_A} F = \{\langle a,b\rangle \in A \times A : (a \leftrightarrow b)^n \in F \textit{ for every } n < \omega\}$.
5. *The logic* $\models_{\mathsf{K}}^{\leqslant}$ *is equivalential.* ∎

Finally, we analyse the issue of the algebraizability of the logics preserving degres of truth. The following fact is interesting to notice: We already know that all logics $\models_{\mathsf{K}}^{\leqslant}$ preserving degrees of truth are selfextensional, and that all logics $\vdash_{\mathsf{K}}$ preserving truth are algebraizable. One reading of the next result is that these properties so-to-speak separate the two groups, in the sense that a logic in one group cannot have the characteristic property of the other group unless it actually belongs to it (because the two logics concide), and that this happens if and only if the class $\mathsf{K}$ is a variety of generalized Heyting algebras.

**Theorem 4.12.** *Let* $\mathsf{K}$ *be a variety of residuated lattices. Then the following conditions are equivalent:*

1. *The logic* $\models_{\mathsf{K}}^{\leqslant}$ *is algebraizable.*
2. $\models_{\mathsf{K}}^{\leqslant} = \vdash_{\mathsf{K}}$.
3. *The logic* $\vdash_{\mathsf{K}}$ *is selfextensional.*
4. $\mathsf{K}$ *is a variety of generalized Heyting algebras (i.e.,* $\mathsf{K} \subseteq \mathsf{Prot_1}$*).*

*Proof.* (1⇒2) Assume $\models_{\mathsf{K}}^{\leqslant}$ is algebraizable and let $\boldsymbol{A} \in \mathsf{K}$; we are going to prove that $\mathcal{F}i_{\models_{\mathsf{K}}^{\leqslant}}(\boldsymbol{A}) = \mathcal{F}i_{\vdash_{\mathsf{K}}}(\boldsymbol{A})$. By (21), in general $\mathcal{F}i_{\models_{\mathsf{K}}^{\leqslant}}(\boldsymbol{A}) \supseteq \mathcal{F}i_{\vdash_{\mathsf{K}}}(\boldsymbol{A})$. Now take any $F \in \mathcal{F}i_{\models_{\mathsf{K}}^{\leqslant}}(\boldsymbol{A})$. In Proposition 3.5 we defined $F^{+}$ as $\min\{G \in \mathcal{F}i_{\models_{\mathsf{K}}^{\leqslant}}(\boldsymbol{A}) : \boldsymbol{\Omega}_{\boldsymbol{A}} G = \boldsymbol{\Omega}_{\boldsymbol{A}} F\}$. Since the Leibniz operator is injective over the filters of an algebraizable logic on an arbitrary algebra, in particular it will be injective over $\mathcal{F}i_{\models_{\mathsf{K}}^{\leqslant}}(\boldsymbol{A})$ and hence $F^{+} = F$. By Proposition 3.5.1 we conclude that $F \in \mathcal{F}i_{\vdash_{\mathsf{K}}}(\boldsymbol{A})$. Since $\mathsf{Alg}(\models_{\mathsf{K}}^{\leqslant}) = \mathsf{Alg}(\vdash_{\mathsf{K}}) = \mathsf{K}$, the equality $\mathcal{F}i_{\models_{\mathsf{K}}^{\leqslant}}(\boldsymbol{A}) = \mathcal{F}i_{\vdash_{\mathsf{K}}}(\boldsymbol{A})$ in all $\boldsymbol{A} \in \mathsf{K}$ implies that $\models_{\mathsf{K}}^{\leqslant} = \vdash_{\mathsf{K}}$.
(2⇒3) because $\models_{\mathsf{K}}^{\leqslant}$ is selfextensional by Proposition 4.1.
(3⇒4) In a residuated lattice $a \star b = 1$ if and only if $a \wedge b = 1$. Therefore, $x \star y \dashv\vdash_{\mathsf{K}} x \wedge y$ for all $\mathsf{K}$. If we assume that $\vdash_{\mathsf{K}}$ is selfextensional, this means that $\mathsf{K} \models x \star y \approx x \wedge y$, and this amounts to saying that $\mathsf{K}$ is a variety of generalized Heyting algebras.
(4⇒2) It is well known that if $\mathsf{K}$ is a variety of generalized Heyting algebras then $\vdash_{\mathsf{K}}$ satisfies the deduction theorem. Applying it to Proposition 2.7 we obtain that $\vdash_{\mathsf{K}}$ and $\models_{\mathsf{K}}^{\leqslant}$ coincide on finite, non-empty sets of assumptions. Since both logics are finitary and by Lemma 2.6 they have the same theorems, they are equal.
(2⇒1) because $\vdash_{\mathsf{K}}$ is always algebraizable. ∎

The reader may have noticed that the only property of algebraizable logics used in the proof of the step 1⇒2 is the injectivity of the Leibniz operator. Thus, this property alone might be listed among the equivalent properties in Theorem 4.12. Moreover, that $\models_{\mathsf{K}}^{\leqslant}$ belongs to any class of logics containing the algebraizable ones and where the Leibniz operator is injective will also be equivalent to any of these properties. The most notable of these, yet not well-known, are the class of *weakly algebraizable* logics [18, 35] and the class of *truth-equational* logics [58]; they are defined, respectively, by the property that the Leibniz operator is both monotone and injective, and by the property that it is completely order-reflecting, a condition which implies injectivity but not monotonicity. It has been proved that algebraizable logics are weakly algebraizable, and that these are truth-equational.

Recall that, by Proposition 4.3, two other properties are also equivalent to those in Theorem 4.12: that the logic $\models_{\mathsf{K}}^{\leqslant}$ is Fregean, and that it is fully Fregean. Since being Fregean implies being selfextensional, from the merging of both results we infer that the algebraizable logic $\vdash_{\mathsf{K}}$ is Fregean if and only if $\mathsf{K}$ is a variety of generalized Heyting algebras. This matches well with Theorem 3.1 in [62], where it is shown that an extension of the logic $FL$ related to the Full Lambek's calculus is Fregean if and only if it is an axiomatic extension of $FL_{eci}$, which is precisely the logic corresponding to the variety of generalized Heyting algebras. It is interesting to notice that the authors of [62] propose a definition of a substructural logic as a non-Fregean extension of $FL$; we do not inted to discuss this proposal, but notice that in such a context our results would characterize the logics $\models_{\mathsf{K}}^{\leqslant}$ for $\mathsf{K}$ a variety of generalized Heyting algebras as non-substructural, or, as we say in the Introduction, as not properly substructural.

The weakest among all logics falling under Theorem 4.12 is (definitionally equivalent to) Johansson's minimal logic. If we add the condition that the constant $0$ is an inconsistent formula (i.e., one that entails all other formulas; this corresponds to all algebras in $\mathsf{K}$ having $0$ as a least element) then we find one of the best known non-classical logics, namely intuitionistic logic. Thus this result yields several characterizations of intuitionistic logic: Among all logics preserving degrees of truth with respect to a variety of bounded residuated lattices, it is the weakest algebraizable one; and among all those preserving truth with respect to a variety of bounded residuated lattices, it is the weakest selfextensional one.

One of the logics that fall under the scope of Theorem 4.12 is the logic corresponding to linear Heyting algebras, known in the literature as either Dummett's logic [23] or as Gödel's logic [45]. Thus in this case $\vdash_{\mathsf{K}}$ and $\models_{\mathsf{K}}^{\leqslant}$ will coincide. This is a strengthening of Theorem 4.2.18 of [45] and of Proposition 13 of [5], which state a similar fact in the case where $\mathsf{K}$ consists of the single algebra on the real unit interval (and, indeed, in the case of [45], its $\models_{\mathsf{K}}^{\leqslant}$ actually uses (1) but only with the rational points, while the model is still the whole real unit interval).

For another family of logics falling under the scope of the previous results, consider the case where $\mathsf{K}$ is the variety generated by a finite subalgebra of the standard real unit interval $[0,1]$; this corresponds to some of the most popular many-valued logics, namely Łukasiewicz's finitely-valued logics. These were analyzed in [42] by means of many-sided sequent calculi, but the associated sentential logics were also considered, and it was proved (with different terminology and notation) that the corresponding logics preserving degrees of truth $\models_{\mathsf{K}}^{\leqslant}$ are finitely equivalential, and non-algebraizable for $n > 2$.

# 5 Tarski and Gentzen style characterizations

In this section we begin by giving (Theorem 5.1) what can be called a Tarski style characterization of the variety $\mathsf{RL}$ of residuated lattices. This term originate with Wójcicki, who in [63] calls *Tarski style condition* any property of a closure operator where only one logical connective or operation appears. By relaxing somehow this definition, i. e., by allowing more than one operation to appear in each condition, and by appying it to closure operators on arbitrary algebras (Wójcicki works only on the set of formulas), we can start by giving a characterization of residuated lattices in terms of Tarski style conditions (notice that actually only one condition is of the "relaxed" kind, namely condition T6 below corresponding to residuation).

The Tarski style conditions will determine a sequent calculus whose reduced models are the operators of lattice filter generation in residuated lattices. From this, we obtain both a Tarski style and a Gentzen style characterization of the logic $\models^{\leqslant}_{\mathsf{RL}}$, which is the weakest case of those treated in a unified form up to now. As we have already noted, its associated substructural logic is the logic denoted as $FL_{ei}$ in [40]. Since all other logics correspond to subvarieties of $\mathsf{RL}$, they are axiomatic extensions of $FL_{ei}$ and admit a similar treatment. At the end of the section we show how this can be done.

Recall from Section 2 that a generalized matrix (g-matrix for short) is a pair $\langle \boldsymbol{A}, \mathcal{C} \rangle$ where $\boldsymbol{A}$ is an algebra and $\mathcal{C}$ is a closure system on its carrier $A$. With each closure system we associate a mapping $C : \mathcal{P}(A) \to \mathcal{P}(A)$, the closure operator defined by $X \subseteq A \ \mapsto\ C(X) = \bigcap \{T \in \mathcal{C} : X \subseteq T\}$. Then $\mathcal{C}$ is the family of $C$-closed subsets of $A$, that is, for each $T \subseteq A$, $T \in \mathcal{C}$ if and only if $C(T) = T$. If $a, a_i \in A$, we write $C(a)$ instead of $C\big(\{a\}\big)$, and $C(a_1, \ldots, a_n)$ instead of $C\big(\{a_1, \ldots, a_n\}\big)$. It is clear that the relation $a \equiv_C b$ defined by $C(a) = C(b)$ is an equivalence relation on $A$.

**Theorem 5.1.** $\boldsymbol{A} \in \mathsf{RL}$ *if and only if there exists a finitary closure operator* $C$ *on* $A$ *satisfying the following properties, for all* $a, b, c \in A$*:*

**T0** *(Reducedness)* $C(a) = C(b)$ *implies* $a = b$.

**T1** *(Identity)* $a \to a \in C(\emptyset)$ *and*

*(Maximum)* $1 \in C(\emptyset)$.

**T2** *(Order) if* $a \to b \in C(\emptyset)$ *then* $b \in C(a)$.

**T3** *(Conjunction)* $C(a \wedge b) = C(a, b)$.

**T4** *(Disjunction)* $C(a \vee b) = C(a) \cap C(b)$.

**T6** *(Residuation)* $c \in C(a \star b)$ *iff* $b \to c \in C(a)$.

**T7** *(Premise permutation)* $C\big(a \to (b \to c)\big) \subseteq C\big(b \to (a \to c)\big)$.

*Proof.* The proof uses a number of technical lemmas that will be proven in Appendix A. If $\boldsymbol{A} \in \mathsf{RL}$ and $C$ is the operator of lattice filter generation, then it is well known (as in every lattice) that it is a finitary closure operator on $A$ satisfying T2, T3 and T4. The other properties reflect some of the fundamental properties of residuated lattices. Conversely, if $C$ is a finitary closure operator satisfying T0–T7, then by T0 the equivalence $\equiv_C$ is the identity, thus the results of the lemmas in Appendix A on $\boldsymbol{A}/{\equiv_C}$ apply directly to $\boldsymbol{A}$. In particular, from item 6 in Lemma A.3 it follows that $\boldsymbol{A} \in \mathsf{RL}$, as was to be proven. ∎

As we will see, property T0 will be treated somehow apart from properties T1–T7: The technical lemmas in Appendix A will not require T0, thus being of a wider applicability.

Some straightforward consequences of some properties among T1–T7, which follow just by using that $C$ is a closure operator, are:

1. $a \in C(b)$ if and only if $C(a) \subseteq C(b)$.
2. The relation $a \equiv_C b$ is an equivalence relation.
3. $C(\emptyset) = C(1) \neq \emptyset$. This follows from T1.

4. $C(\emptyset) = \bigcap_{a\in A} C(a)$. This is a consequence of 3.
5. Any finitely generated closed set is principal: $C(a_1,\dots,a_n) = C(a_1\wedge\cdots\wedge a_n)$ is a generalized form of of T3.
6. $C\big(a\to(b\to c)\big) = C\big(b\to(a\to c)\big)$. This follows from T7, by symmetry.

It is worth highlighting here some of the results obtained in Appendix A. For instance, in Lemma A.3 we see:

**Corollary 5.2.** *If a closure operator $C$ on an arbitrary algebra $\boldsymbol{A}$ satisfies properties* T1–T7*, then the equivalence relation $\equiv_C$ is a congruence of $\boldsymbol{A}$ and the quotient $\boldsymbol{A}/{\equiv_C}$ is a residuated lattice.* ∎

**Corollary 5.3.** *If $\boldsymbol{A}$ is a residuated lattice, then the only finitary closure operator $C$ satisfying properties* T0–T7 *in Theorem 5.1 is the one corresponding to the closure system $\mathcal{C} = \mathcal{F}i_{\wedge}(\boldsymbol{A})$.*

*Proof.* It is well known that the closure operation of lattice filter generation in a residuated lattice satisfies the stated properties. Conversely, property T0 means that the relation $\equiv_C$ is the identity, so that $\boldsymbol{A}$ can be identified with $\boldsymbol{A}/{\equiv_C}$; then from Lemma A.4 it follows that if $C$ is finitary and satisfies T1–T7, then $\mathcal{C} = \mathcal{F}i_{\wedge}(\boldsymbol{A})$. ∎

Next we introduce a ***Gentzen system*** for the variety $\mathsf{RL}$. Here a ***sequent*** will be a pair $\langle \Gamma, \varphi\rangle$ where $\Gamma$ is a finite, possibly empty *set* of formulas, and $\varphi$ a formula; we will denote it by $\Gamma \rhd \varphi$ in order to avoid any potential misunderstanding with other symbols that are sometimes used as sequent separator, such as $\vdash$, $\to$ or $\Rightarrow$.

**Definition 5.4.** *The Gentzen system $\mathcal{G}_{\mathsf{RL}}$ is the one defined by the following axioms and rules (where $\varphi, \psi, \xi$ range over arbitrary formulas):*

**1.** *All the structural axioms and rules, i.e., the initial axiom $\varphi \rhd \varphi$ plus the rules of weakening and cut.*

**2.** *Logical axioms:*

| | |
|---|---|
| $\emptyset \rhd \varphi\to\varphi$ | *(Identity)* |
| $\emptyset \rhd 1$ | *(Maximum)* |
| $\varphi\,,\psi \rhd \varphi\wedge\psi$ | *(Conjunction 1)* |
| $\varphi\wedge\psi \rhd \varphi$ | *(Conjunction 2)* |
| $\varphi\wedge\psi \rhd \psi$ | *(Conjunction 3)* |
| $\psi\to(\varphi\to\xi) \rhd \varphi\to(\psi\to\xi)$ | *(Premise permutation)* |

**3.** *Logical rules:*

$$\frac{\emptyset \rhd \varphi\to\psi}{\varphi \rhd \psi}\ \textit{(Order)} \qquad \frac{\varphi\vee\psi \rhd \xi}{\varphi \rhd \xi}\ \textit{(Disjunction 1)}$$

$$\frac{\varphi\vee\psi \rhd \xi}{\psi \rhd \xi}\ \textit{(Disjunction 2)} \qquad \frac{\varphi \rhd \xi \quad \psi \rhd \xi}{\varphi\vee\psi \rhd \xi}\ \textit{(Disjunction 3)}$$

$$\frac{\varphi\star\psi \rhd \xi}{\varphi \rhd \psi\to\xi}\ \textit{(Residuation 1)} \qquad \frac{\varphi \rhd \psi\to\xi}{\varphi\star\psi \rhd \xi}\ \textit{(Residuation 2)}$$

Since the left-hand side of our sequents is a (finite) set of formulas, rather than a multiset or a sequence, it is not necessary to include the structural rules of contraction and exchange, as they are so-to-speak built-in in the formalism.

**Definition 5.5.** *A generalized matrix $\langle \boldsymbol{A}, \mathcal{C} \rangle$ with corresponding closure operator $C$ is a* ***model of a Gentzen system*** *when it is a model of all its rules, in the following sense: Let*

$$\frac{\Gamma_0 \rhd \varphi_0 \quad \ldots \quad \Gamma_{n-1} \rhd \varphi_{n-1}}{\Gamma_n \rhd \varphi_n} \tag{27}$$

*be the general form of a Gentzen style rule. $\langle \boldsymbol{A}, \mathcal{C} \rangle$ is a* ***model of the rule*** (27) *when for any $v \in \mathrm{Hom}(\boldsymbol{Fm}, \boldsymbol{A})$, if $v(\varphi_i) \in C\big(v[\Gamma_i]\big)$ for all $i < n$, then $v(\varphi_n) \in C\big(v[\Gamma_n]\big)$.*

Obviously, axioms (initial sequents) are treated as rules with no premises, in which case the above definition yields that $\langle \boldsymbol{A}, \mathcal{C} \rangle$ is a model of an axiom $\Gamma \rhd \psi$ when $v(\psi) \in C\big(v[\Gamma]\big)$ for all $v \in \mathrm{Hom}(\boldsymbol{Fm}, \boldsymbol{A})$. Notice that any closure operator is a model of the structural rules. If $L$ is a logic then we will say that $L$ **satisfies the rule** (27) when $\langle \boldsymbol{Fm}, \vdash_L \rangle$ is a model of the rule.

The rules of our Gentzen system have been chosen (and labelled) so as to make the proof of the following result simply trivial:

**Theorem 5.6.** *Let $\boldsymbol{A}$ be an algebra, and $C$ a closure operator on $\boldsymbol{A}$. Then the associated g-matrix $\langle \boldsymbol{A}, \mathcal{C} \rangle$ is a model of $\mathcal{G}_{\mathsf{RL}}$ if and only if $C$ satisfies the properties* T1–T7 *of Theorem 5.1.* ∎

**Proposition 5.7.** *If $\boldsymbol{A}$ is an algebra and $C$ is a finitary closure operator on $\boldsymbol{A}$, then the associated g-matrix $\langle \boldsymbol{A}, \mathcal{C} \rangle$ is a reduced model of $\mathcal{G}_{\mathsf{RL}}$ if and only if $\boldsymbol{A} \in \mathsf{RL}$ and $\mathcal{C} = \mathcal{F}i_{\wedge}(\boldsymbol{A})$.*

*Proof.* If $\boldsymbol{A} \in \mathsf{RL}$ and $\mathcal{C} = \mathcal{F}i_{\wedge}(\boldsymbol{A})$ then Corollary 5.3 plus Theorem 5.6 imply that $\langle \boldsymbol{A}, C \rangle$ is a model of $\mathcal{G}_{\mathsf{RL}}$, and moreover it is reduced by condition T0. Conversely, if $\langle \boldsymbol{A}, C \rangle$ is a model of $\mathcal{G}_{\mathsf{RL}}$ then by Theorem 5.6 $C$ satisfies properties T1–T7. Then by Corollary 5.2 $\equiv_C$ is a congruence. Now the assumption of being reduced means that $\equiv_C$ is the identity, thus $C$ satisfies T0, and then Corollary 5.3 finishes the proof. ∎

Since the algebraic counterpart of a Gentzen system is defined, as for deductive systems, as the class of algebra reducts of its reduced models, we have:

**Corollary 5.8.** $\mathsf{Alg}(\mathcal{G}_{\mathsf{RL}}) = \mathsf{RL}$ *and hence also* $\mathsf{Alg}(\mathcal{G}_{\mathsf{RL}}) = \mathsf{Alg}(\models_{\mathsf{RL}}^{\leqslant})$. ∎

From the preceding results a very compact description of the relation between the Gentzen system $\mathcal{G}_{\mathsf{RL}}$ and the logic $\models_{\mathsf{RL}}^{\leqslant}$ can be obtained. The notion of a Gentzen system being ***fully adequate*** for a sentential logic, introduced in [35, Definition 4.10], has a deep significance in some recent works in abstract algebraic logic, see [37, 39]. When the logic has theorems (which is the case here), then it can be defined by saying that the full g-models of the logic coincide with the models of the Genten system. It is important to stress that this means that the Gentzen system characterizes not only all the rules satisfied by the logic, but also those inherited by the operator of filter generation in arbitrary algebras.

**Theorem 5.9.** *The Gentzen system* $\mathcal{G}_{\mathsf{RL}}$ *is fully adequate for* $\models^{\leqslant}_{\mathsf{RL}}$.

*Proof.* In Corollary 5.8 we have proved that $\mathcal{G}_{\mathsf{RL}}$ and $\models^{\leqslant}_{\mathsf{RL}}$ have the same algebraic counterpart. Since in $\boldsymbol{A} \in \mathsf{RL}$, by Proposition 2.9, $\mathcal{F}i_{\wedge}(\boldsymbol{A}) = \mathcal{F}i_{\models^{\leqslant}_{\mathsf{RL}}}(\boldsymbol{A})$, Corollary 5.3 and Theorem 5.1 put together show that on residuated lattices the only reduced model of $\mathcal{G}_{\mathsf{RL}}$ corresponds to the closure system of all filters of $\models^{\leqslant}_{\mathsf{RL}}$. This is exactly the characterization obtained in Proposition 4.12 of [35] of the notion of a Gentzen system being fully adequate for a deductive system, in the case where the latter has theorems, which is the case here. ∎

Other consequences that follow immediately from Theorem 5.9 are:

**Corollary 5.10.**

1. $\models^{\leqslant}_{\mathsf{RL}}$ *is the internal deductive system associated with* $\mathcal{G}_{\mathsf{RL}}$ *in the following sense: if* $\{\varphi_0, \dots, \varphi_n\} \subseteq Fm$, *then*

$$\varphi_0, \dots, \varphi_{n-1} \models^{\leqslant}_{\mathsf{RL}} \varphi_n \iff \varphi_0, \dots, \varphi_{n-1} \rhd \varphi_n \text{ is derivable in } \mathcal{G}_{\mathsf{RL}}.$$

2. $\models^{\leqslant}_{\mathsf{RL}}$ *is the weakest deductive system that satisfies all the axioms and rules of* $\mathcal{G}_{\mathsf{RL}}$.
3. $\models^{\leqslant}_{\mathsf{RL}}$ *is the weakest deductive system that satisfies properties* T1–T7 *of Theorem 5.1.*
4. $\mathcal{G}_{\mathsf{RL}}$ *is algebraizable with respect to* $\mathsf{RL}$ *with transformers*

$$\begin{aligned} \emptyset \rhd \psi &\longmapsto \psi \approx 1 \\ \varphi_0, \dots, \varphi_{n-1} \rhd \varphi_n &\longmapsto \varphi_0 \wedge \cdots \wedge \varphi_{n-1} \wedge \varphi_n \approx \varphi_0 \wedge \cdots \wedge \varphi_{n-1} \\ \varphi \approx \psi &\longmapsto \{\varphi \rhd \psi, \psi \rhd \varphi\}. \end{aligned}$$

5. $\vdash_{\mathsf{RL}}$ *is the external deductive system associated with* $\mathcal{G}_{\mathsf{RL}}$ *in the following sense: if* $\{\varphi_0, \dots, \varphi_n\} \subseteq Fm$, *then* $\varphi_0, \dots, \varphi_{n-1} \vdash_{\mathsf{RL}} \varphi_n$ *if and only if the sequent* $\emptyset \rhd \varphi_n$ *is derivable in* $\mathcal{G}_{\mathsf{RL}}$ *from the sequents* $\emptyset \rhd \varphi_0, \dots, \emptyset \rhd \varphi_{n-1}$. ∎

The terms "internal" and "external" are taken from [4]. Part 2 is what can be called a *Gentzen style characterization* of $\models^{\leqslant}_{\mathsf{RL}}$, while part 3 is the *(almost) Tarski style characterization* of $\models^{\leqslant}_{\mathsf{RL}}$. Part 5 can be paraphrased as saying that $\vdash_{\mathsf{RL}}$ is *the Hilbert subrelation* of the Gentzen system $\mathcal{G}_{\mathsf{RL}}$ in the sense of [57]; moreover, since both the logic $\vdash_{\mathsf{RL}}$ and the Gentzen system $\mathcal{G}_{\mathsf{RL}}$ are algebraizable and the variety $\mathsf{RL}$ is their common equivalent algebraic semantics, it follows that $\vdash_{\mathsf{RL}}$ and $\mathcal{G}_{\mathsf{RL}}$ are *deductively equivalent*, and that $\mathcal{G}_{\mathsf{RL}}$ is *simply Hilbertizable*, again using the terminology of [57]. It is interesting to notice the twofold relationship between the variety $\mathsf{RL}$ and the logical systems analysed in this paper. While by definition the variety $\mathsf{RL}$ gives raise to the two sentential logics $\vdash_{\mathsf{RL}}$ and $\models^{\leqslant}_{\mathsf{RL}}$ in a semantic way, it also defines the Gentzen system $\mathcal{G}_{\mathsf{RL}}$, which in turn also happens to determine syntactically the two logics as shown above.

There is another Gentzen system having a further interest for $\vdash_{\mathsf{RL}}$, namely the substructural system originating in [54] and algebraically studied in [1], which is a variant of the one called $FL_{ew}$ in the literature; it is shown that this Gentzen system is also algebraizable, having $\mathsf{FL_{ew}}$ as its equivalent algebraic semantics, and is also deductively equivalent to $\vdash_{\mathsf{FL_{ew}}}$ and to $\mathcal{G}_{\mathsf{FL_{ew}}}$. By eliminating the axiom concerning the constant 0 we obtain its analogues for $\mathsf{RL}$, with parallel properties.

Finally, let $\mathsf{K}$ be a subvariety of $\mathsf{RL}$. In some cases, an equational presentation of $\mathsf{K}$ relatively to $\mathsf{RL}$ is known. Besides the cases already mentioned in Section 2, other examples are all the varieties of BL-algebras generated by a finite number of continuous t-norms[4], which together with other varieties of BL-algebras have been axiomatized in [21, 26, 49]; all the varieties of nilpotent minimum algebras have been axiomatized in [43]. In all these cases, all the results found in this section for $\models^{\leqslant}_{\mathsf{RL}}$ can be extended to $\models^{\leqslant}_{\mathsf{K}}$. It suffices to add the condition $C(\varphi) = C(\psi)$ for each equation $\varphi \approx \psi$ of the said presentation to the conditions in Theorem 5.1, and we obtain a Tarski style characterization of the variety $\mathsf{K}$. And by adding the two sequential axioms $\varphi \rhd \psi$ and $\psi \rhd \varphi$ to Definition 5.4 we obtain a Gentzen system $\mathcal{G}_{\mathsf{K}}$ for which all the remaining results hold, *mutatis mutandis*. Thus, we obtain a Tarski style and a Gentzen style characterizations of $\models^{\leqslant}_{\mathsf{K}}$. For instance, to treat the cases where 0 is the minimum of the algebra we should add the condition "$C(0) = A$" to Theorem 5.1, and the logical axiom "$0 \rhd \varphi$" to Definition 5.4. Let us mention that this is not the only possible way of carrying out this investigation, and for particular cases other *ad hoc* presentations can be found; this was for instance the case of Łukasiewicz's infinite-valued logic as treated in [33].

If no equational presentation of $\mathsf{K}$ over $\mathsf{RL}$ is known, one can still do the same for each equation $\varphi \approx \psi$ that holds in $\mathsf{K}$ and not in $\mathsf{RL}$, but this will then yield an infinite, and perhaps non-recursive, family of conditions and axioms.

# A Proof of some lemmas

In this appendix we give the proofs of the technical lemmas necessary for the proof of Theorem 5.1 and its corollaries. These lemmas do not need condition T0. The first one concerns only the operations of implication and fusion.

**Lemma A.1.** *If $C$ satisfies properties* T1, T2, T6, T7 *then for any $a, b, c \in A$:*

1. *if $\{a, a \to b\} \subseteq C(\emptyset)$ then $b \in C(\emptyset)$, i.e., $C(\emptyset)$ is closed under Modus Ponens.*
2. $a \in C(b \star a)$.
3. *if $c \in C(\emptyset)$ then $a \to c \in C(\emptyset)$.*
4. *if $C(a) \subseteq C(b)$ then $C(a \star c) \subseteq C(b \star c)$ and $C(c \to a) \subseteq C(c \to b)$.*
5. *$b \in C(a)$ iff $a \to b \in C(\emptyset)$.*
6. *$C(a) = C(b)$ iff $\{a \to b, b \to a\} \subseteq C(\emptyset)$.*
7. $C(a \star b) = C(b \star a)$.

[4] That is, algebras whose universe is the real interval $[0, 1]$, the lattice operations are those corresponding to the natural order, $\star$ is given by a (continuous) t-norm, and $\to$ by its residual.

**8.** *if* $C(a) \subseteq C(b)$ *then* $C(b \to c) \subseteq C(a \to c)$.

**9.** $a \to (b \to c) \in C(\emptyset)$ *iff* $(a \star b) \to c \in C(\emptyset)$.

**10.** $C((a \star b) \to c) = C(a \to (b \to c))$.

**11.** $C(a \star (b \star c)) = C((a \star b) \star c)$.

**12.** $C((a \to a) \to b) = C(b)$.

*Proof.* 1. By T2 and the hypotheses, $a \to b \in C(\emptyset)$ implies $b \in C(a) \subseteq C(\emptyset)$.

2. By T1 and T6, $a \to a \in C(\emptyset) \subseteq C(b)$ implies $a \in C(b \star a)$.
3. If $c \in C(\emptyset) \subseteq C(c \star a)$, then by T6 $a \to c \in C(c) = C(\emptyset)$.
4. By T6, $d \in C(a \star c)$, and so $c \to d \in C(a) \subseteq C(b)$; hence $d \in C(b \star c)$. So $C(a \star c) \subseteq C(b \star c)$. For the other inclusion, note that by T6 it holds that $c \to b \in C(c \to b)$. Hence $a \in C(b) \subseteq C((c \to b) \star c)$. Therefore, $c \to a \in C(c \to b)$.
5. ($\Leftarrow$) is T2. ($\Rightarrow$) If $a \in C(b)$ then $C(a) \subseteq C(b)$, so by 4 $C(b \to a) \subseteq C(b \to b) = C(\emptyset)$.
6. It is a consequence of 5.
7. Since $a \star b \in C(a \star b)$, by T6 and 5, $a \to (b \to (a \star b)) \in C(\emptyset)$. Now by T7 $b \to (a \to (a \star b)) \in C(\emptyset)$. Finally, by 5 and T6 again, $a \star b \in C(b \star a)$. Symmetrically, $b \star a \in C(a \star b)$ and hence $C(a \star b) = C(b \star a)$.
8. The proof is similar to the one given in the second part of 4. By T6, 7, 4 and T6 again we have that $a \to c \in C(a \to c)$. Thus, $c \in C((a \to c) \star a) \subseteq C((a \to c) \star b)$, and hence $b \to c \in C(a \to c)$, so $C(b \to c) \subseteq C(a \to c)$.
9. By 5, T6 and 5 again, $a \to (b \to c) \in C(\emptyset)$ iff $b \to c \in C(a)$ iff $c \in C(a \star b)$ iff $(a \star b) \to c \in C(\emptyset)$.
10. Note that

$$
\begin{array}{rll}
C((a \to (b \to c)) \to ((a \star b) \to c)) \in C(\emptyset) & \text{iff} & \text{by T7} \\
C((a \star b) \to ((a \to (b \to c)) \to c)) \in C(\emptyset) & \text{iff} & \text{by 9} \\
C(a \to (b \to ((a \to (b \to c)) \to c))) \in C(\emptyset) & \text{iff} & \text{by T7 and 4} \\
C(a \to ((a \to (b \to c)) \to (b \to c))) \in C(\emptyset) & \text{iff} & \text{by T7} \\
C((a \to (b \to c)) \to (a \to (b \to c))) \in C(\emptyset) & & \text{by T1,}
\end{array}
$$

and that

$$
\begin{array}{rll}
((a \star b) \to c) \to (a \to (b \to c)) \in C(\emptyset) & \text{iff} & \text{by T7} \\
a \to (((a \star b) \to c) \to (b \to c)) \in C(\emptyset) & \text{iff} & \text{by 5} \\
((a \star b) \to c) \to (b \to c) \in C(a) & \text{iff} & \text{by T7} \\
b \to (((a \star b) \to c) \to c) \in C(a) & \text{iff} & \text{by 5} \\
a \to (b \to (((a \star b) \to c) \to c)) \in C(\emptyset) & \text{iff} & \text{by 10} \\
(a \star b) \to (((a \star b) \to c) \to c)) \in C(\emptyset) & \text{iff} & \text{by T7} \\
((a \star b) \to c) \to ((a \star b) \to c)) \in C(\emptyset) & & \text{by T1.}
\end{array}
$$

Applying 6 to the two previously proven facts, we obtain that $C(a \to (b \to c)) = C((a \star b) \to c)$.

11. Note that

$$
\begin{aligned}
C(((a\star b)\star c)\to(a\star(b\star c))) &= &&\text{by 10}\\
C(a\star b\to(c\to(a\star(b\star c)))) &= &&\text{by 10}\\
C(a\to(b\to(c\to(a\star(b\star c))))) &= &&\text{by 10 and 4}\\
C(a\to((b\star c)\to(a\star(b\star c)))) &= &&\text{by 10}\\
C((a\star(b\star c))\to(a\star(b\star c))) &= C(\emptyset) &&\text{by T1.}
\end{aligned}
$$

so by T2 $a\star(b\star c)\in C((a\star b)\star c))$; to prove $(a\star b)\star c)\in C((a\star(b\star c))$ we proceed in a similar way.

12. By T1, T6, 7 and T6, we have $b\to b\in C(a\to a)$, so $b\in C((a\to a)\star b)=C(b\star(a\to a))$, and then $(a\to a)\to b\in C(b)$ by T6. The proof of $b\in C(a\to a)\to b)$ is similar. ∎

Remark that property 5 establishes the converse of the implication in T2, hence the equivalence. Thus, the relation "$a\to b\in C(\emptyset)$" is a quasi-order. This relation can also be expressed in terms of the lattice-like operations:

**Lemma A.2.** *If $C$ satisfies* T1–T7 *then for all $a,b\in A$,*

$$a\to b\in C(\emptyset)\ \textit{ iff }\ C(a)=C(a\wedge b)\ \textit{ iff }\ C(b)=C(a\vee b).$$

*Proof.* By 5 and T3, $a\to b\in C(\emptyset)$ iff $b\in C(a)$ iff $C(a)=C(a,b)=C(a\wedge b)$. For the second equivalence, note that by 5 and T4, $a\to b\in C(\emptyset)$ iff $b\in C(a)$ iff $C(b)\subseteq C(a)$ iff $C(b)=C(a)\cap C(b)=C(a\vee b)$. ∎

Recall that the relation $\equiv_C$ is always an equivalence relation. Now we prove that with our properties it is a congruence relation, and see what happens in the quotient:

**Lemma A.3.** *If $C$ satisfies* T1–T7 *then:*

**1.** *$a\equiv_C b$ if and only if $a\to b\in C(\emptyset)$ and $b\to a\in C(\emptyset)$.*

**2.** *$\equiv_C$ has the substitution property relative to $\to$ and $\star$.*

**3.** *For any $a\in A$, $(a\to a)/{\equiv_C}=1/{\equiv_C}=C(\emptyset)$.*

**4.** *The relation $a/{\equiv_C}\leqslant b/{\equiv_C}\iff C(b)\subseteq C(a)\iff a/{\equiv_C}\to b/{\equiv_C}=1/{\equiv_C}$ is a partial order in $A/{\equiv_C}$, and $1/{\equiv_C}$ is its maximum, which is denoted equally by 1.*

**5.** *$\equiv_C\in \mathrm{Co}\boldsymbol{A}$ (it is a congruence of the algebra $\boldsymbol{A}$).*

**6.** *$\boldsymbol{A}/{\equiv_C}$ is a residuated lattice.*

*Proof.* 1. This is another expression of property 6 in Lemma A.1.

2. This follows from properties 4, 7 and 8 in Lemma A.1.
3. This is a consequence of T1 and the (general) fact that all elements in $C(\emptyset)$ are equivalent.
4. The first equivalence and that this relation is a partial order in $A/{\equiv_C}$ are also general facts of the theory of closure operators. The second equivalence follows from property 5 in Lemma A.1, given that by T1, $C(\emptyset)=C(1)$. That $1/{\equiv_C}$ is the maximum of $A/{\equiv_C}$ follows from property 3 in Lemma A.1.

5. The substitution property relative to $\to$ and $\star$ has already been established, and it is well known that the substitution properties relative to $\wedge$ and $\vee$ are a consequence of T3 and T4, respectively.
6. Properties T3 and T4 imply that $\boldsymbol{A}/{\equiv_C}$ is a lattice, and we have just established (part 4) that 1 is its maximum. Property T6 says, modulo part 4, that $\star$ is residuated and $\to$ is its residual. Concerning $\star$, properties 11 and 7 in Lemma A.1, modulo part 5, show that it is a commutative semigroup. Finally, to show that it is a monoid we have to show that 1 is the unit, and by part 5 it is enough to show that $C(a \star 1) = C(a)$: By properties 2 and 7 in Lemma A.1, $C(a) \subseteq C(a \star 1)$. On the other hand, by T1, $(a \star 1) \to (a \star 1) \in C(\emptyset)$, so by parts 7 and 9 of Lemma A.1, $1 \to (a \to (a \star 1)) \in C(\emptyset)$, and this, after using property 5 of the same lemma twice together with the fact that $C(1) = C(\emptyset)$ by T1, gives $a \star 1 \in C(a)$, that is, $C(a \star 1) \subseteq C(a)$. ∎

Given $T \subseteq A$, we define $T/{\equiv_C} = \{a/{\equiv_C} : a \in T\}$. Then:

**Lemma A.4.** *If $C$ satisfies* T1–T7 *then:*

**1.** *If $T \in \mathcal{C}$ then $T = \bigcup T/{\equiv_C}$.*

**2.** $C(a)/{\equiv_C} = \{b/{\equiv_C} : a/{\equiv_C} \leqslant b/{\equiv_C}\}$.

**3.** *If moreover $C$ is finitary, then the map $T \mapsto T/{\equiv_C}$ is an order isomorphism from the family of closed sets of $C$ onto the family of lattice filters of $A/{\equiv_C}$, both ordered by inclusion.*

*Proof.* 1. In general, if $a \in T$ then $a \in \bigcup T/{\equiv_C}$ because $a \in a/{\equiv_C}$. Conversely, assume that $T \in \mathcal{C}$, that is, $T = C(T)$. If $a \in \bigcup T/{\equiv_C}$, this means that $a \in b/{\equiv_C}$ for some $b \in T$, so $C(b) \subseteq T$ and $a \equiv_C b$, that is, $C(a) = C(b)$ and then $a \in T$.

2. This follows from property 4 in Lemma A.3.

3. If $T \in \mathcal{C}$ then from T3 it follows that $T/{\equiv_C}$ is a lattice filter. Moreover, by 1, $T/{\equiv_C} \subseteq T'/{\equiv_C}$ if and only if $T \subseteq T'$. This also implies that the mapping is one-to-one. Now let $F$ be a lattice filter of $\boldsymbol{A}/{\equiv_C}$ and put $T_F = C\big(\bigcup F\big)$. Then, if $a \in T_F$ there exist $a_1, \dots, a_n \in \bigcup F$ such that $a \in C(a_1, \dots, a_n) = C(a_1 \wedge \dots \wedge a_n)$, so $a/{\equiv_C} \geqslant (a_1 \wedge \dots \wedge a_n)/{\equiv_C} = a_1/{\equiv_C} \wedge \dots \wedge a_n/{\equiv_C} \in F$, which implies that $a \in \bigcup F$ because $F$ is a lattice filter. Thus $T_F = \bigcup F$. This shows that the map is onto and completes the proof. ∎

# B Some examples

In this appendix we present five residuated lattices that have been used in the paper as counterexamples to several statements. The underlying lattice structures are indicated by the Hasse diagrams in Figure 1, while the remaining operations (i.e., fusion and implication) are introduced in each case.

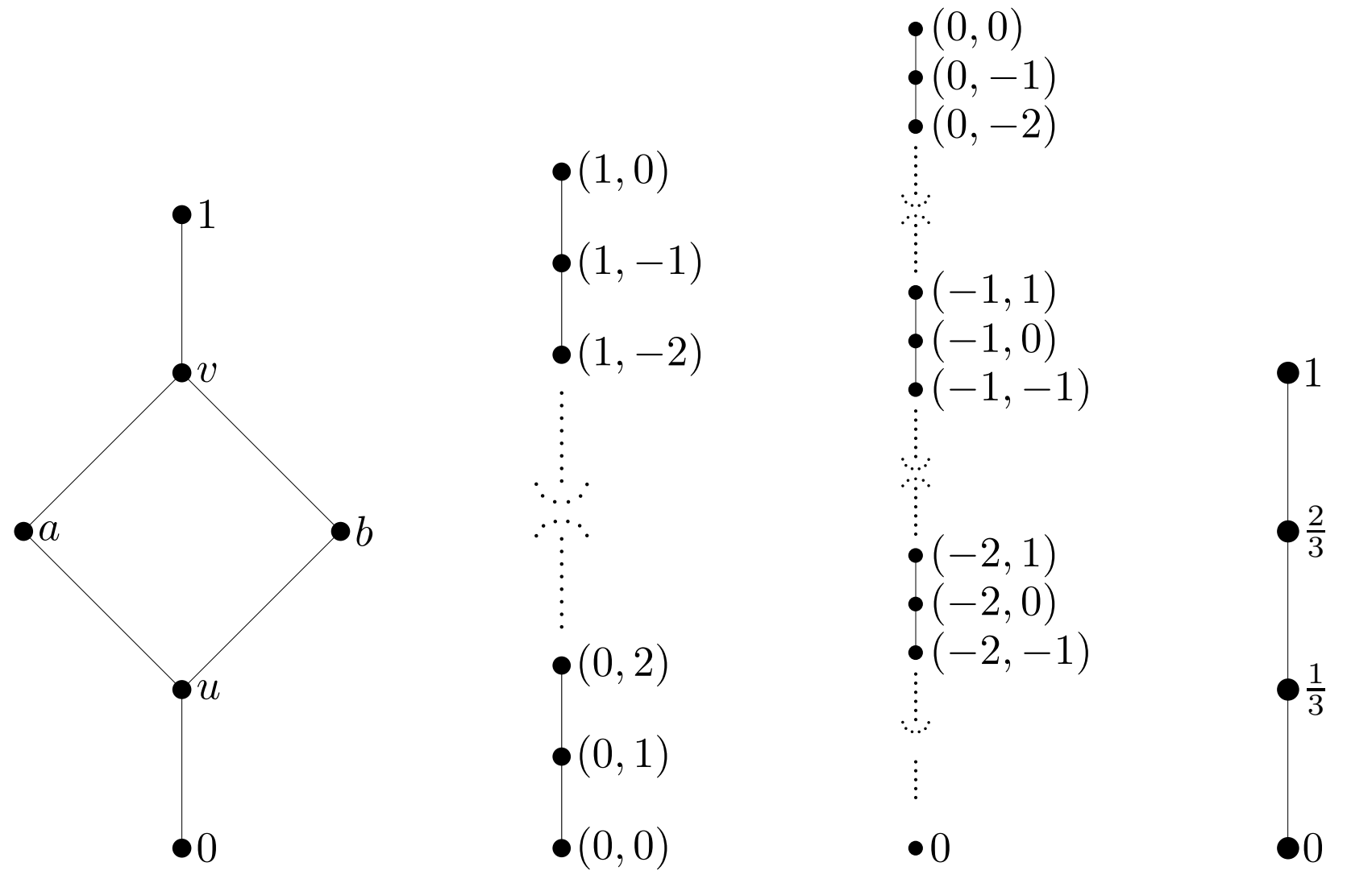

(a) $\boldsymbol{A}_1$ and $\boldsymbol{A}_4$ (b) $\boldsymbol{A}_2$ (c) $\boldsymbol{A}_3$ (d) $\boldsymbol{A}_5$

Figure 1: Underlying lattices of the examples in Appendix B.

**Example 1**

The residuated lattice $\boldsymbol{A_1}$ has a 6-element universe $\{0, u, a, b, v, 1\}$ and the lattice structure depicted in Figure 1(a). The fusion and implication operations are given by the following tables:

| $\star$ | $0$ | $u$ | $a$ | $b$ | $v$ | $1$ |
|---|---|---|---|---|---|---|
| $0$ | $0$ | $0$ | $0$ | $0$ | $0$ | $0$ |
| $u$ | $0$ | $0$ | $0$ | $0$ | $0$ | $u$ |
| $a$ | $0$ | $0$ | $a$ | $0$ | $a$ | $a$ |
| $b$ | $0$ | $0$ | $0$ | $b$ | $b$ | $b$ |
| $v$ | $0$ | $0$ | $a$ | $b$ | $v$ | $v$ |
| $1$ | $0$ | $u$ | $a$ | $b$ | $v$ | $1$ |

| $\to$ | $0$ | $u$ | $a$ | $b$ | $v$ | $1$ |
|---|---|---|---|---|---|---|
| $0$ | $1$ | $1$ | $1$ | $1$ | $1$ | $1$ |
| $u$ | $v$ | $1$ | $1$ | $1$ | $1$ | $1$ |
| $a$ | $b$ | $b$ | $1$ | $b$ | $1$ | $1$ |
| $b$ | $a$ | $a$ | $a$ | $1$ | $1$ | $1$ |
| $v$ | $u$ | $u$ | $a$ | $b$ | $1$ | $1$ |
| $1$ | $0$ | $u$ | $a$ | $b$ | $v$ | $1$ |

This algebra is used in the paper for two different purposes. First, because the lattice filter generated by $u$ (i.e., $\{u, a, b, v, 1\}$) does not contain a maximum implicative filter (cf. Theorem 4.4). Second, because $\boldsymbol{A_1}$ witnesses that $\mathsf{IMC_n}$ is strictly included into $\mathsf{E_n}$ (Theorem 4.9, part 2), since $a$ and $b$ are idempotent elements while $a \wedge b$ is not (cf. Lemma 4.7).

**Example 2**

This residuated lattice corresponds to the MV-algebra known as "Chang's algebra" in the literature, see [11, pag. 481] for instance. The universe of the algebra $\boldsymbol{A_2}$ is

$$\big\{(0, x) : x \in \mathbb{Z}^+\big\} \cup \big\{(1, x) : x \in \mathbb{Z}^-\big\}.$$

where $\mathbb{Z}^+$ denotes the set of positive integers, including $0$, while $\mathbb{Z}^-$ denotes the negative ones, also including $0$. These elements are ordered as shown in Figure 1(b). The fusion $\star$ is defined by

$$(i,x) \star (j,y) := \begin{cases} (0,0) & \text{if } i=j=0 \\ (1, x+y) & \text{if } i=j=1 \\ (0, \max\{0, x+y\}) & \text{otherwise} \end{cases}$$

and the implication $\to$ is defined by

$$(i,x) \to (j,y) := \begin{cases} (1,1) & \text{if } i=0 \text{ and } j=1 \\ (0, -x+y) & \text{if } i=1 \text{ and } j=0 \\ (1, \min\{0, -x+y\}) & \text{otherwise} \end{cases}$$

As an MV-algebra $\boldsymbol{A_2}$ is the one obtained through Mundici's "Gamma" functor from the lattice-ordered abelian group $\mathbb{Z} \times_\ell \mathbb{Z}$, where $\times_\ell$ is lexicographical product, taking $(1,0)$ as unit, see [12, Section 2.1]. Chang's algebra gives us an example of an MV-algebra where the Leibniz operator is not monotone over lattice filters, i.e., condition (M) fails, thus showing that the logic $\models^{\leqslant}_{\mathsf{MV}}$ is non-protoalgebraic. To see this, it is enough to consider the lattice filters $F = \big\{(1,x) : x \in \mathbb{Z}^-\big\}$ and $G = F \cup \big\{(0,x) : x \geqslant 1\big\}$. Since $F$ is an implicative filter, we know that $\boldsymbol{\Omega}_{\boldsymbol{A_2}} F$ is the congruence associated with $F$, and an easy examination shows that $\boldsymbol{\Omega}_{\boldsymbol{A_2}} G$ is the congruence associated with the implicative filter $\big\{(1,0)\big\}$. Therefore, $\big\langle (0,0)\,,(0,1)\big\rangle \in \boldsymbol{\Omega}_{\boldsymbol{A_2}} F$ while $\big\langle (0,0)\,,(0,1)\big\rangle \notin \boldsymbol{\Omega}_{\boldsymbol{A_2}} G$. Thus, $F \subseteq G$ but $\boldsymbol{\Omega}_{\boldsymbol{A_2}} F \nsubseteq \boldsymbol{\Omega}_{\boldsymbol{A_2}} G$.

**Example 3**

The universe of the algebra $\boldsymbol{A_3}$ is

$$\big\{(0,x) : x \in \mathbb{Z}^-\big\} \cup \big\{(y,x) : y \in \mathbb{Z}^- \smallsetminus \{0\}, x \in \mathbb{Z}\big\} \cup \big\{0\big\},$$

and its lattice structure is given in Figure 1(c). Its fusion $\star$ is defined by

$$(i,x) \star (j,y) := (i+j\,, x+y) \quad \text{and} \quad z \star 0 = 0 \star z = 0$$

while its implication $\to$ is given by

$$\begin{aligned} (i,x) \to 0 \quad &:= 0 \\ 0 \to (j,y) \;&:= (0,0) \\ (i,x) \to (j,y) \;&:= \begin{cases} \big(0, \min\{0, y-x\}\big) & \text{if } i \leqslant j \\ (j-i\,, y-x) & \text{if } j < i \end{cases} \end{aligned}$$

This is a product algebra, namely the one associated in [13] with the negative cone of the linearly ordered abelian group $\mathbb{Z} \times_\ell \mathbb{Z}$. The algebra $\boldsymbol{A_3}$ is an example of a product algebra where the Leibniz operator is not monotone over lattice filters. The reasoning is similar to that in Example 2, now using the implicative filter $F = \big\{(0,x) : x \in \mathbb{Z}^-\big\}$ and the lattice filter $G = F \cup \big\{(-1,x) : x \geqslant 0\big\}$. It is not difficult to check that $\big\langle (-1,0)\,,(-1,1)\big\rangle \in \boldsymbol{\Omega}_{\boldsymbol{A_3}} F$ while $\big\langle (-1,0)\,,(-1,1)\big\rangle \notin \boldsymbol{\Omega}_{\boldsymbol{A_3}} G$.

**Example 4**

The residuated lattice $\boldsymbol{A_4}$ has the same universe and lattice structure as that of Example 1, shown in Figure 1(a). The fusion and implication are given by the following tables:

| $\star$ | 0 | $u$ | $a$ | $b$ | $v$ | 1 |
|---|---|---|---|---|---|---|
| 0 | 0 | 0 | 0 | 0 | 0 | 0 |
| $u$ | 0 | $u$ | $u$ | $u$ | $u$ | $u$ |
| $a$ | 0 | $u$ | $a$ | $u$ | $a$ | $a$ |
| $b$ | 0 | $u$ | $u$ | $u$ | $u$ | $b$ |
| $v$ | 0 | $u$ | $a$ | $u$ | $a$ | $v$ |
| 1 | 0 | $u$ | $a$ | $b$ | $v$ | 1 |

| $\to$ | 0 | $u$ | $a$ | $b$ | $v$ | 1 |
|---|---|---|---|---|---|---|
| 0 | 1 | 1 | 1 | 1 | 1 | 1 |
| $u$ | 0 | 1 | 1 | 1 | 1 | 1 |
| $a$ | 0 | $b$ | 1 | $b$ | 1 | 1 |
| $b$ | 0 | $v$ | $v$ | 1 | 1 | 1 |
| $v$ | 0 | $b$ | $v$ | $b$ | 1 | 1 |
| 1 | 0 | $u$ | $a$ | $b$ | $v$ | 1 |

This algebra is an example of a residuated lattice (more precisely, a $\mathrm{FL}_{ew}$ algebra) where the set of idempotent elements, here $\{0, u, a, 1\}$, is closed under meet, while it is not an MTL algebra; this is so because $(a \to b) \vee (b \to a) \neq 1$. Moreover, $\boldsymbol{A_4} \in \mathsf{Prot_2}$, but $\boldsymbol{A_4}$ is not an ordinal sum of simple n-contractive algebras; this can be seen by using that in simple algebras there are no non-trivial idempotent elements, analysing all possibilities.

**Example 5**

The universe of $\boldsymbol{A_5}$ is the 4-element set $\left\{0\,, \frac{1}{3}\,, \frac{2}{3}\,, 1\right\}$, with the lattice structure given by the natural linear ordering, as shown in Figure 1(d). The fusion and implication in $\boldsymbol{A_5}$ are the operations given by

| $\star$ | 0 | $\frac{1}{3}$ | $\frac{2}{3}$ | 1 |
|---|---|---|---|---|
| 0 | 0 | 0 | 0 | 0 |
| $\frac{1}{3}$ | 0 | 0 | 0 | $\frac{1}{3}$ |
| $\frac{2}{3}$ | 0 | 0 | $\frac{2}{3}$ | $\frac{2}{3}$ |
| 1 | 0 | $\frac{1}{3}$ | $\frac{2}{3}$ | 1 |

| $\to$ | 0 | $\frac{1}{3}$ | $\frac{2}{3}$ | 1 |
|---|---|---|---|---|
| 0 | 1 | 1 | 1 | 1 |
| $\frac{1}{3}$ | $\frac{2}{3}$ | 1 | 1 | 1 |
| $\frac{2}{3}$ | $\frac{1}{3}$ | $\frac{1}{3}$ | 1 | 1 |
| 1 | 0 | $\frac{1}{3}$ | $\frac{2}{3}$ | 1 |

This algebra is a subalgebra of the algebra given by the nilpotent minimum t-norm (see [25, 43]). It holds that $\boldsymbol{A_5} \in \mathsf{MTL} \cap \mathsf{IMC_2} \subseteq \mathsf{MTL} \cap \mathsf{IMC_n}$ while $\boldsymbol{A_5} \notin \mathsf{Prot_n}$ because $\frac{1}{3} \wedge \left((\frac{1}{3} \to 0)^n \star (0 \to \frac{1}{3})^n\right) \nleq 0$.

FÉLIX BOU
Artificial Intelligence Research Institute (IIIA - CSIC)
Bellaterra
fbou@iiia.csic.es

FRANCESC ESTEVA
Artificial Intelligence Research Institute (IIIA - CSIC)
Bellaterra
esteva@iiia.csic.es

JOSEP MARIA FONT
Faculty of Mathematics
University of Barcelona
jmfont@ub.edu

ÀNGEL J. GIL
Departament d'Economia i Empresa
Universitat Pompeu Fabra, Barcelona
angel.gil@upf.edu

LLUÍS GODO
Artificial Intelligence Research Institute (IIIA - CSIC)
Bellaterra
godo@iiia.csic.es

ANTONI TORRENS
Faculty of Mathematics
University of Barcelona
atorrens@ub.edu

VENTURA VERDÚ
Faculty of Mathematics
University of Barcelona
v.verdu@ub.edu